\documentclass[12pt, leqno]{amsart}

\usepackage[OT2,T1]{fontenc}
\usepackage{indentfirst}
\usepackage{amstext}
\usepackage{amsthm}
\usepackage{amsopn}
\usepackage{amsfonts}
\usepackage{amsmath}
\usepackage{latexsym}
\usepackage[francais,english]{babel}
\usepackage{amscd}
\usepackage{amssymb}
\usepackage{amsmath}
\usepackage[all,cmtip]{xy}
\usepackage{leftidx}
\usepackage{graphicx}
\usepackage{etoolbox}
\patchcmd{\section}{\normalfont\scshape\centering}{\normalfont\bfseries}{}{}
\patchcmd{\subsection}{-.5em}{.5em}{}{}

\newtheorem{theo}{{Theorem}}[section]
\newtheorem{coro}[theo]{{Corollary}}
\newtheorem{lemma}[theo]{{Lemma}}
\newtheorem{prop}[theo]{Proposition}

\theoremstyle{definition}
\newtheorem{remark}[theo]{\textbf{Remark}}
\newtheorem{defn}[theo]{Definition}
\newtheorem{example}[theo]{Example}
\numberwithin{equation}{section}

\newtheorem{notation}[theo]{Notation}

\DeclareFontEncoding{OT2}{}{} 
  \newcommand{\textcyr}[1]{%
    {\fontencoding{OT2}\fontfamily{wncyr}\fontseries{m}\fontshape{n}%
     \selectfont #1}}
\newcommand{\sha}{{\mbox{\textcyr{Sh}}}}

\begin{document}
\tolerance 400 \pretolerance 200 \selectlanguage{english}

\title{Isometries of lattices and automorphisms of $K3$ surfaces}
\author{Eva  Bayer-Fluckiger}
\date{\today}
\maketitle

\begin{abstract} The aim of this paper is to give necessary and sufficient conditions for an integral polynomial  to be the characteristic
polynomial of a semi-simple isometry of some even unimodular lattice of given signature. This result has applications to
automorphisms of 
$K3$ surfaces; in particular, we show that every Salem number of degree $4,6,8,12,14$ or $16$ is the dynamical
degree of an automorphism of a non-projective $K3$ surface. 

\medskip


\medskip

\end{abstract}

\small{} \normalsize

\medskip

\selectlanguage{english}
\section{Introduction}

Let $r,s \geqslant 0$ be integers such that $r  \equiv s \ {\rm (mod \ 8)}$; this congruence condition is equivalent to the
existence of an even, unimodular lattice with signature $(r,s)$. When $r,s \geqslant 1$, such a lattice is unique up
to isomorphism (see for instance \cite{S}, chap. V); we denote it by $\Lambda_{r,s}$.
In \cite {GM}, Gross and McMullen raise
the following question (see \cite {GM}, Question 1.1) :

\bigskip
\noindent
{\bf Question.} {\it What are the possibilities for the characteristic polynomial $F(X) = {\rm det}(X - t)$ of an isometry
$t \in {\rm SO}(\Lambda_{r,s})$ ? }

\bigskip
The aim of this paper is to answer this question for semi-simple isometries.

\bigskip
The condition $t \in {\rm SO}(\Lambda_{r,s})$ implies that $F(X) = X^{{\rm deg}(F)}F(X^{-1})$,
hence $F$ is  a {\it symmetric} polynomial (cf. \S \ref{symmetric section}). Let
$2n = {\rm deg}(F)$, and let $2m(F)$ be the number of roots of $F$ outside the unit circle. 
As shown in \cite{GM}, we have the further necessary conditions :

\bigskip
(C 1) {\it $|F(1)|$, $|F(-1)|$ and $(-1)^n F(1) F(-1)$ are squares.}

\bigskip

(C 2) {\it $r \geqslant m(F)$, $s \geqslant m(F)$, and if moreover $F(1)F(-1) \not = 0$, then $m(F) \equiv r \equiv s  \ {\rm (mod \ 2)}$.}

\bigskip
Gross and McMullen prove that if $F \in {\bf Z}[X]$ is an irreducible, symmetric and monic
polynomial satisfying condition (C 2) and such that  $|F(1)F(-1)| = 1$, then there exists  $t \in {\rm SO}(\Lambda_{r,s})$ with characteristic polynomial $F$ (see
\cite {GM}, Theorem 1.2). They
speculate that conditions (C 1) and (C 2) are sufficient for a monic {\it irreducible} polynomial to be realized as the characteristic
polynomial of an isometry of $\Lambda_{r,s}$; this is proved in \cite {BT}, Theorem A. More generally, Theorem A of \cite {BT}
implies that 
if a monic, irreducible 
and symmetric polynomial $F$ satisfies conditions (C 1) and (C 2), then {\it there exists} an even, unimodular lattice
of signature $(r,s)$ having an isometry with characteristic polynomial $F$. This is also the point of view of the present
paper - we treat the definite and indefinite cases simultaneously. 

\medskip
On the other hand, Gross and McMullen show that these conditions do not suffice in the case
of {\it reducible} polynomials (see \cite {GM}, Proposition 5.2); several other examples are given in \cite{B 20}. Another
example is the following :

\medskip
\noindent
{\bf Example 1.} Let $F(X) = (X^4-X^2+1)(X-1)^4$, and let $(r,s) = (8,0)$; conditions (C 1) and (C 2) hold, but there does not exist any positive definite, even, unimodular lattice of rank $8$ having an isometry with characteristic polynomial $F$; note that
this amounts to saying that the lattice $E_8$ does not have any isometry with characteristic polynomial $F$. 

\medskip
All these examples are {\it counter-examples to a Hasse principle}. Indeed, 
the first result of the present paper is that conditions (C 1) and (C 2) are sufficient {\it locally}. 
 If $p$ is a prime number, we say
that a ${\bf Z}_p$-lattice $(L,q)$ is even if $q(x,x) \in 2{\bf Z}_p$ for all $x \in L$; note that if $p \not = 2$, then every lattice is even, since $2$
is a unit in ${\bf Z}_p$. The following is proved in Theorem \ref{integral local} and  Proposition \ref{reals}~:

\medskip
\noindent
{\bf Theorem 1.} {\it Let $F \in {\bf Z}[X]$ be a monic, symmetric polynomial of even degree. 

\medskip  {\rm (a)} Condition {\rm (C 1)} holds if and only if for all prime numbers $p$, there
exists an even,  unimodular ${\bf Z}_p$-lattice having a semi-simple isometry
with characteristic polynomial $F$.

\medskip {\rm (b)} 
The group ${\rm SO}_{r,s}({\bf R})$ contains a semi-simple element having characteristic polynomial $F$  if and
only if condition {\rm (C 2)} holds.}

\medskip
The next result is a necessary and sufficient condition for the local-global principle to hold. We start
by defining an obstruction group (see \S \ref{obstruction section}). Let us write $F(X) = F_1(X)(X-1)^{n_+}(X+1)^{n_-}$ and
assume that $n_+ \not = 2$, $n_- \not = 2$; in this case, the group only depends on the polynomial $F$; we denote it by
$\sha_F$.

 \medskip Let us now assume that condition (C 2) holds, and let $t \in {\rm SO}_{r,s}({\bf R})$ be a semi-simple
 isometry with characteristic polynomial $F$; such an isometry exists by part (b) of the above theorem (or Proposition 
 \ref{reals}). Assume moreover that condition (C 1) also holds. In \S \ref{reformulation},  we define a homomorphism 
 $$\epsilon_t : \sha_F \to {\bf Z}/2{\bf Z}$$ and prove the following (see Theorem \ref{preserves}) : 
 
 \medskip
 \noindent
 {\bf Theorem 2.} {\it The isometry  $t \in {\rm SO}_{r,s}({\bf R})$ preserves an even, unimodular lattice if and
 only if $\epsilon_t = 0$.}
 
 \bigskip
 \noindent
 {\bf Example 2.} Let $F(X) = (X^4-X^2+1)(X-1)^4$; we have $\sha_F \simeq {\bf Z}/2{\bf Z}$.
 If $t_1 \in {\rm SO}_{4,4}({\bf R})$, we have $\epsilon_{t_1} = 0$, and if $t_2 \in {\rm SO}_{8,0}({\bf R})$, then
 $\epsilon_{t_2} \not = 0$.  Hence $\Lambda_{4,4}$ has a semi-simple isometry with characteristic polynomial $F$,
 but the lattice $E_8$ does not have such an isometry.

 \bigskip
 \noindent
 {\bf Corollary 1.} {\it  
 Let $G \in {\bf Z}[X]$ be a monic,  irreducible,  symmetric polynomial such that $|G(1)|$ is not a square, and suppose that $|G(-1)|$ is a square. Let
 $m \geqslant 4$  be an even integer, and set $$F(X) = G(X)(X-1)^{m}.$$ 
 
 \medskip
 \noindent
 Assume that condition {\rm (C 2)} holds for $F$. 
 Then every semi-simple isometry   $t \in {\rm SO}_{r,s}({\bf R})$ with characteristic polynomial $F$ preserves an even, unimodular lattice.}
 
 \medskip Indeed, Condition (C 1) holds for $F$ since $|G(-1)|$ is a square, and one can check that $\sha_F = 0$;  therefore Theorem 2 implies the corollary.

 \bigskip
 For polynomials $F$ without linear factors, Theorem 2 is proved in \cite{B 20}, Theorem 27.4. However, it turns out that including
 linear factors is very useful in the applications to $K3$ surfaces, which we now describe. 
 
 \medskip The second part of the paper gives applications to automorphisms of $K3$ surfaces,
 inspired by a series of papers of McMullen  (see \cite {Mc1}, \cite {Mc2}, \cite {Mc3}).
  
 \medskip
 
 Recall that a
 monic, irreducible, symmetric polynomial $S \in {\bf Z}[X]$ of degree $\geqslant 4$ is a {\it Salem polynomial} if $S$ has exactly two roots outside
 the unit circle, both positive real numbers. A real number is called a {\it Salem number} if it is the unique real root $> 1$ of
 a Salem polynomial; it is an algebraic unit.

 \medskip If 
 $T : \mathcal X \to \mathcal X$ is an automorphism of a complex $K3$ surface, then 
 $T^* : H^2(\mathcal X,{\bf C}) \to H^2(\mathcal X,{\bf C})$ respects the Hodge decomposition
 $$H^2(\mathcal X,{\bf C}) = H^{2,0}(\mathcal X) \oplus H^{1,1}(\mathcal X) \oplus H^{0,2}(\mathcal X);$$ since
 ${\rm dim}(H^{2,0}) = 1$, the automorphism $T^*$ acts on it by multiplication with a complex number, denoted by
 $\delta(T)$; we have $|\delta(T)| = 1$. Moreover, $T^* : H^2(\mathcal X,{\bf Z}) \to H^2(\mathcal X,{\bf Z})$ preserves
 the intersection pairing.
 The above properties imply that the characteristic polynomial
of $T^*$ is a product of at most one Salem polynomial and of a finite number of cyclotomic polynomials, 
it satisfies condition (C 1), and $\delta(T)$ is a root of this polynomial. 

\medskip
Moreover, assume that the characteristic polynomial is equal to $S C$, where $S$ is a Salem polynomial
of degree $d$ with $4 \leqslant d \leqslant 22$ and $C$ is a product of cyclotomic polynomials; then $\mathcal X$
is projective if and only if $\delta(T)$ is a root of $C$ (see \cite{R2}, Theorem 2.2). Such a polynomial is called a {\it complemented
Salem polynomial} (see Definition \ref{complemented}). 

\medskip 

  Let $F$ be a complemented Salem polynomial, 
  and let 
 $\delta$ be a root of $F$. We wish
  to decide whether $F$ is the characteristic polynomial of an isomorphism $T^*$ as above, with $\delta(T) = \delta$. 
   We start with
  the simplest case, where the cyclotomic factor is a nontrivial power of $X-1$.

 \bigskip
 \noindent
 {\bf Theorem 3.}
 {\it Let $S$ be a Salem polynomial of degree $d$ with $4 \leqslant d\leqslant 18$
 and let $\delta$ be a root of $S$ with $|\delta| = 1$. Let $$F(X) = S(X)(X-1)^{22-d},$$
 assume that condition {\rm (C 1)} holds for $F$, and that $|S(1)|$ is not a square.
 Then
 there exists a non-projective $K3$ surface $\mathcal X$
 and an automorphism $T : \mathcal X \to \mathcal X$ such
 that
 
 \medskip \noindent
 $\bullet$
 \ \ $F$ is the characteristic polynomial of $T^*|H^2(\mathcal X).$ 
 
  \medskip \noindent
$\bullet$ \ \ $T^*$ acts on $H^{2,0}(\mathcal X)$ by multiplication by $\delta$.}

\medskip
This is proved using Corollary 1, as well as some results of McMullen (\cite {Mc2}, {Mc3}) and Brandhorst (\cite {Br}).
For polynomials $S$ with $|S(1)| =1$, see Theorem \ref{theorem 4 - second part}; in
this case, the answer depends on the congruence class of $d$ modulo $8$.

\medskip The {\it dynamical degree} of an automorphism $T : \mathcal X \to \mathcal X$ is by definition the spectral
radius of $T^*$; since the characteristic polynomial of $T^*$ is the product of a Salem polynomial and of a product
of cyclotomic polynomial, the dynamical degree is a Salem number.  We say that a Salem number is 
{\it realizable} if $\alpha$ is the dynamical degree of an automorphism of a 
 $K3$ surface. 

\medskip
Let $\alpha$ be a Salem number of degree $d$ with $4 \leqslant d\leqslant 20$, and let $S$ be the minimal polynomial of
$\alpha$. In \S \ref{4} we prove
an analog of Theorem 3 for $F(X) = S(X)(X+1)^2(X-1)^{20-d}$ or $S(X)(X^2-1)(X-1)^{20-d}$, and show 
that if $d = 4,6,8,12,14$ or $16$, then $\alpha$ is realizable (see Corollary \ref{lots coro}).

\medskip

The aim of \S \ref{18}  is to prove that the second smallest known Salem number, $\lambda_{18} = 1.1883681475...$, is not
realizable as the dynamical degree of an automorphism of a  {\it non-projective} $K3$ surface.  
By contrast, McMullen proved that $\lambda_{18} $ is the dynamical degree
of an automorphism of a {\it projective} $K3$ surface (see \cite {Mc3}, Theorem 8.1).

\bigskip
I thank Marie Jos\'e Bertin, Serge Cantat, Curt McMullen, Chris Smyth and Yuta Takada for very useful comments and suggestions.

\medskip

\section{Equivariant Witt groups}

\medskip
We start by recalling some notions and results from \cite {BT}, \S 3 and \S 4. 

\medskip
{\bf The equivariant Witt group}

\medskip 
Let $\mathcal K$ be a field, let $\mathcal A$ be a $\mathcal K$-algebra and let $\sigma : \mathcal A \to \mathcal A$ be a $\mathcal K$-linear involution. 
An {\it $\mathcal A$-bilinear form} is a pair $(V,b)$ consisting of an $\mathcal A$-module $V$ that is a finite dimensional $\mathcal K$-vector space,
and a non-degenerate symmetric $\mathcal K$-bilinear form $b : V \times V \to K$ such that $b(ax,y) = b(x,\sigma(a)y)$ for all $a \in \mathcal A$ and
all $x, y \in V$. 

\medskip
The associated {\it Witt group} is denoted by $W_{\mathcal A}(\mathcal K)$ (see \cite {BT}, \S 3). If $M$ is a simple $\mathcal A$-module, we
denote by $W_{\mathcal A}({\mathcal K},M)$ the subgroup of $W_{\mathcal A}(\mathcal K)$ generated by the classes of $\mathcal A$-bilinear
forms $(M,b)$. Every class in $W_{\mathcal A}(\mathcal K)$ is represented by an $\mathcal A$-bilinear form whose underlining $\mathcal A$-module
is semisimple, and we have $$W_{\mathcal A}(\mathcal K) = \underset{M} \oplus W_{\mathcal A}({\mathcal K},M),$$ where $M$ ranges over the
isomorphism classes of simple $\mathcal A$-modules (see \cite {BT}, Corollary 3.11 and Theorem 3.12).

\medskip
{\bf Discrete valuation rings and residue maps}

\medskip
Let $O$ be a discrete valuation ring with field of fractions $K$, residue field and uniformizer $\pi$. Let $(A,\sigma)$ be an $O$-algebra with involution,
and set $A_K = A \otimes_O K$, $A_k = A \otimes_O k$. An {\it $A$-lattice} in an $A_K$ bilinear form $V$ is an $A$-submodule $L$ which
is finitely generated as an $O$-module and satisfies $K L = V$. If $L$ is an $A$-lattice, then so is its dual $$L^{\sharp} = \{ x \in V \ | \ b(x,L) \subset O \}.$$
We say that $L$ is {\it unimodular} if $L^{\sharp} = L$ and {\it almost unimodular} if $\pi L^{\sharp} \subset L \subset L^{\sharp}$. If $L$ is almost unimodular,
then $b$ induces an $A_k$-bilinear form $L^{\sharp}/L \times L^{\sharp}/L \to k$ (see \cite {BT}, definition 4.1).

\medskip

An $A_K$-bilinear form
is said to be {\it bounded} if it contains an $A$-lattice. We denote by $W^b_{A_K}(K)$ the subgroup of $W_{A_K}(K)$ generated by the classes of
bounded forms. The following result is proved in \cite{BT} :

\begin{theo}\label{4.3}  
{\rm (i)} Every bounded $A_K$-bilinear form contains an almost unimodular $A$-lattice $L$.

\medskip
{\rm (ii)} The class of $L^{\sharp}/L$ in $W_{A_k}(k)$ only depends on the class of $V$ in $W_A(K)$. 

\medskip
{\rm (iii)} The map $\partial : W_{A_K}(K) \to W_{A_k}(k)$ given by $[V] \to [L^{\sharp}/L]$ is a homomorphism.

\medskip
{\rm (iv)} $V$ contains a unimodular $A$-lattice if and only if $V$ is bounded and $\partial [V] = 0$ in $W_{A_k}(k)$.

\end{theo}

\noindent
{\bf Proof.} See \cite {BT}, Theorem 4.3. 

\section{Symmetric polynomials and $\Gamma$-modules}\label{symmetric section}

We recall some notions from \cite {M} and  \cite {B 15}. Let $K$ be a field. If
$f \in K[X]$ is a monic polynomial such that $f(0) \not = 0$,  set $f^*(X) = {f(0)} X^{{\rm deg}(f)} f(X^{-1})$; we say that $f$ is
{\it symmetric} if $f^* = f$.  Recall the following definition from \cite {B 15} :

\begin{defn} Let $f \in k[X]$ be a monic, symmetric polynomial. We say that $f$ is of 

\medskip
\noindent
$\bullet$ {\it type} 0 if $f$ is a product of powers of $X-1$ and of $X+1$;

\medskip
\noindent
$\bullet$ {\it type} 1  if $f$ is a product of powers of monic, symmetric, irreducible polynomials in
$k[X]$ of even degree;

\medskip
\noindent
$\bullet$ {\it type} 2 if $f$ is a product of polynomials of the form $g g^*$, where $g \in k[X]$ is
monic, irreducible, and $g \not = \pm  g^*$. 

\end{defn}

The following is well-known :

\begin{prop} Every monic symmetric polynomial  is a product of polynomials
of type {\rm 0, 1} and {\rm 2}.

\end{prop}

\noindent
{\bf Proof.} See for instance  \cite {B 15}, Proposition 1.3.

\medskip 
Let $J$ be the set of irreducible factors of $F$, and let us write $F = \underset{f \in J} \prod f^{n_f}$. Let $I_1 \subset J$ be the subset of irreducible factors of type 1, and let $I_0 \subset J$ be the
set of irreducible factors of type 0; set $I = I_0 \cup I_1$. For all $f \in I$, set $M_f = [ K[X]/(f)]^{n_f}$. Set $M^0 =  \underset {f \in I_{0}} \oplus M_f$, and
$M^1 =  \underset {f \in I_1} \oplus M_f$. 
If $f \in J$ such that $f \not = f^*$, set $M_{f,f^*} = [K[X]/(f) \oplus K[X]/(f^*)]^{n_f}$, and let $M^2 = \underset {(f,f^*)} \oplus M_{f,f^*}$,
where the sum runs over the pairs $(f,f^*)$ with $f \in J$ and $f \not = f^*$.  Set $$M = M^0 \oplus M^1 \oplus M^2.$$

\medskip
Let $\Gamma$ be the infinite cyclic group, and let $\gamma$ be a generator of $\Gamma$. Setting $\gamma (m) = X m$ for all $m \in M$ endows $M$
with a structure of semi-simple $K[\Gamma]$-module; we say that $M$ is the {\it semi-simple $K[\Gamma]$-module associated to the polynomial
$F$}. 

\medskip Let us write $F = F_0 F_1 F_2$, where $F_i$ is the product of the irreducible factors of type $i$ of $F$. We have
$F_0 = (X-1)^{n^+} (X+1)^{n^-}$ for some integers $n^+,n^- \geqslant 0$. Set $M^+ = [K[X]/(X-1)^{n^+}$ and 
$M^- = [K[X]/(X+1)^{n^-}$. The $K[\Gamma]$-module $M^0$ splits as $$M^0 = M^+ \oplus M^-.$$

\section{Isometries of quadratic forms}\label{isometries section}

We recall some results from \cite {M} and \cite {B 15}. Let $K$ be a field of characteristic $\not = 2$, let $V$ be a finite dimensional
$K$-vector space, and let $q : V \times V \to K$ be a non-degenerate quadratic form. An {\it isometry} of $(V,q)$ is by definition
an isomorphism $t : V \to V$ such that $q(tx,ty) = q(x,y)$ for all $x, y \in V$. 
Let $t : V \to V$ be an isometry, and let $F \in K[X]$ be the characteristic polynomial of $t$. It is well-known that $F$ 
is a symmetric polynomial (see for instance \cite {B 15}, Proposition 1.1). The following property is also well-known :

\begin{lemma}\label{determinant} If $t : V \to V$ is an isometry of the quadratic form $(V,q)$ and if the characteristic
polynomial $F$ of $t$ satisfies $F(1)F(-1) \not = 0$, then $${\rm det}(q) = F(1)F(-1) \ \ {\rm in} \ \ K^{\times}/K^{\times 2}.$$

\end{lemma}

\noindent
{\bf Proof.} See for instance \cite{B 15}, Corollary 5.2.

\medskip
Recall that $\Gamma$ is the infinite cyclic group,
and let $\sigma : K[\Gamma] \to K[\Gamma]$ be the $K$-linear involution such that
$\sigma(\gamma) = 
\gamma^{-1}$ for all $\gamma \in \Gamma$. An isometry $t  : V \to V$ endows $V$ with a $K[\Gamma]$-module structure,
and if moreover $t$ is semi-simple with characteristic polynomial $F$, then 
this module is isomorphic to the
semi-simple $K[\Gamma]$-module $M = M(F)$ associated to the polynomial
$F$ (see \S \ref{symmetric section}). Hence $M$ also carries a non-degenerate quadratic form, that we also denote by $q$.
Note that $(M,q)$ is a $K[\Gamma]$-bilinear form, and gives rise to an element $[M,q]$ of
the Witt group $W_{K[\Gamma]}(K)$. To simplify notation, set $W_{\Gamma}(K) = W_{K[\Gamma]}(K)$. 

\medskip
Let us write $M = M^0 \oplus M^1 \oplus M^2$ as in \S \ref{symmetric section}, and let $q^i$ denote the restriction of $q$ to $M^i$; this gives
rise to an orthogonal decomposition $(M,q) = (M^0,q^0) \oplus (M^1,q^1) \oplus (M^2,q^2)$, and $(M^2,q^2)$ is hyperbolic, hence its class
in $W_{\Gamma}(K)$ is trivial (see for instance \cite {M}, Lemma 3.1). With the notation of \S \ref{symmetric section}, we have the further
orthogonal decompositions 

\medskip 

\centerline {$(M^0,q^0) = \underset {f \in I_0} \oplus (M_f,q_f)$ and  $(M^1,q^1) = \underset {f \in I_1} \oplus (M_f,q_f)$,}

\noindent
where $q_f$ is the restriction of $q$ to $M_f$ (see for instance \cite {M}, \S 3, or \cite {B 15},
 Propositions 3.3 and 3.4).  Note that if $f \in I_0$, then
$f(X) = X-1$ or $X+1$, and we have the orthogonal decomposition $(M^0,q^0) = (M^+,q^+) \oplus (M^-,q^-)$, with
$q^+ = q_{X-1}$ and $q^- = q_{X+1}$.


\section{Local fields and unimodular $\Gamma$-lattices}\label{local}

\medskip

Let $K$ be a non-archimedean local field of characteristic 0, let $O$ be its ring of integers, and let $k$ be its residue field.
If  $a \in O$, set $v(a) = 1$ if $v_K(a)$ is  odd, and $v(a) = 0$ if $v_K(a)$ is even or $a = 0$ (in other words, $v(a)$ is the valuation of $a$ (mod 2) if
$a \not = 0$, and $v(0) = 0$). 

\medskip
\begin{theo} \label{local odd}  Let $F \in O[X]$ be a monic, symmetric polynomial. There exists  a unimodular $O$-lattice having a
semi-simple isometry with characteristic polynomial $F$ if and only if one of the following holds 

\medskip {\rm (i)} ${\rm char}(k) \not = 2$, and
 $v(F(1)) = v(F(-1)) = 0$. 
 
\medskip {\rm (ii)} ${\rm char}(k) = 2$, and
 $v(F(1) F(-1)) = 0$.

\end{theo}

\medskip We start with a preliminary result, and some notation. 

\begin{notation}\label{lambda}  Let $E_0$ be an \'etale $K$-algebra of finite rank, and let $E$ be an \'etale $E_0$-algebra 
which is free of rank 2 over $E_0$. Let $\sigma : E \to E$ be the involution fixing $E_0$. If $\lambda \in E_0^{\times}$, we denote by $b_{\lambda}$
the quadratic form $b_{\lambda} : E \times E \to K$ such that $b_{\lambda}(x,y) = {\rm Tr}_{E/K}(\lambda x \sigma(y))$. 

\end{notation}

\begin{prop} \label{technical}
Let $E_0$ be an \'etale $K$-algebra of finite rank, and let $E$ be an \'etale $E_0$-algebra 
which is free of rank 2 over $E_0$. Let $\sigma : E \to E$ be the involution fixing $E_0$.  Let $\alpha \in E_0^{\times}$ be such that
$\alpha \sigma(\alpha) = 1$, and that the characteristic polynomial $f$ of $\alpha$ over $K$ belongs to $O[X]$. Let ${\rm deg}(f) = 2d$, and assume that $f(1) f(-1) \not = 0$. Let
 $u_+, u_-  \in O^{\times}$. 

\medskip
Let $V = V^{+} \oplus V^{-}$ be a finite dimensional $K$-vector space, and let 
$\epsilon = (\epsilon^+,\epsilon^-) : V \to V$ be the isomorphism given by $\epsilon^{\pm} : V^{\pm} \to V^{\pm}$, $\epsilon^{\pm} = \pm id$. 
Set $n^+ = {\rm dim}(V^+)$ and $n^- = {\rm dim}(V^-)$.

\medskip
If ${\rm char}(k) \not = 2$, assume that if $n^{\pm} = 0$, then $v(f(\pm 1)) = 0$.

\medskip
If ${\rm char}(k) = 2$, assume that if $n^+ = n^- = 0$, then $v(f(1)) + v(f(-1)) = 0$. 

\medskip
If moreover $K = {\bf Q}_2$, assume that 

\medskip
$\bullet$
$n^+$  and $n^-$ are both even,   

\medskip
$\bullet$  if $n^+ = n^- = 0$, then
$(-1)^df(1)f(-1) = 1$ in  ${\bf Q}_2^{\times}/{\bf Q_2}^{\times 2}$, 

\medskip
$\bullet$ if $n^{\pm} = 0$, then $v(f(\pm 1)) = 0$,

\medskip
$\bullet$ $u_+ u_- = (-1)^n$, where $2n = {\rm dim}(E \oplus V)$. 


\bigskip
Then there exists $\lambda \in E_0^{\times}$ and non-degenerate quadratic forms 
$$q^+ : V^+ \times V^+ \to K, \ \ q^- : V^- \times V^- \to K$$ such that, for $q = q^+ \oplus q^-$, we have 

\medskip

{\rm (i)}   $$\partial [E\oplus V, b_{\lambda} \oplus q, \alpha \oplus \epsilon] = 0$$
in $W_{\Gamma}(k)$.

\medskip

{\rm (ii)} If  ${\rm char}(k) \not = 2$ 
then ${\rm  det}(q^{\pm}) = u_{\pm}f(\pm 1)$ in $K^{\times}/K^{\times 2}$.

\medskip

{\rm (iii)} If moreover $K = {\bf Q}_2$,  then  

\medskip
$\bullet$ If $n_- \not = 0$, then $v({\rm det}(q^-)) = v(f(-1))$,

\medskip
$\bullet$  If $n_+ \not = 0$ and $n_- \not = 0$, then ${\rm det}(q_{\pm}) = u_{\pm} f(\pm 1)$ in ${\bf Q}_2^{\times}/{\bf Q_2}^{\times 2}$,

\medskip
$\bullet$ ${\rm det}(E\oplus V, b \oplus q) = (-1)^n$, 
and $(E\oplus V, b \oplus q)$ contains an even, unimodular ${\bf Z}_2$-lattice.

\end{prop}

\medskip

\noindent{\bf Proof.} The proof depends on the values of $v(f(1))$ and $v(f(-1))$. We
are in one of the following cases

\medskip
{\rm (a)} $v(f(1)) = 0$, $v(f(-1)) = 0$,

\medskip
{\rm (b)} $v(f(1)) = 1$, $v(f(-1)) = 0$,

\medskip
(c) $v(f(1)) = 0$, $v(f(-1)) = 1$,

\medskip
(d) $v(f(1)) = 1$, $v(f(-1)) = 1$.

\medskip

The algebra $E_0$ decomposes as
a product of fields $E_0 = \underset {v \in S} \prod E_{0,v}$. For all $v \in S$, set $E_v = E \otimes_{E_0} E_{0,v}$. 

\medskip Assume first that the characteristic of $k$ is $ \not = 2$. The algebra $E_v$ is of one of
the following types

\medskip
(sp) $E_v = E_{0,v} \times E_{0,v}$;
 
\medskip
(un) $E_v$ is an unramified extension of $E_{0,v}$;

\medskip
(+) $E_v$ is a ramified extension of $E_{0,v}$, and the image $\overline \alpha$ of $\alpha$ in the residue field $\kappa_v$ of $E_v$ is $1$;

\medskip
(-) $E_v$ is a ramified extension of $E_{0,v}$, and the image $\overline \alpha$ of $\alpha$ in the residue field $\kappa_v$ of $E_v$ is $-1$.

\medskip
This gives a partition $S = S_{sp} \cup S_{un} \cup S_+ \cup S_-$.

\medskip
Let $\gamma$ be a generator of $\Gamma$, and let $\chi_{\pm} : \Gamma \to \{ \pm \}$ be the character sending $\gamma$ to $\pm 1$. 

\medskip 
Let us choose $\lambda = (\lambda_v)_{v \in S}$ in $E_0^{\times} = \underset {v \in S} \prod E_{0,v}^{\times}$ such that 
for every $v \in S_{un}$, we have $\partial [E_v,b_{\lambda_v},\alpha] = 0$ in $W_{\Gamma}(k$); this is possible by \cite {BT}, Proposition 6.4. 
The choices for $v \in S_+$ and $S_-$ depend on which of the cases (a), (b), (c) or (d) we are in. Let $\overline u$ be the image of $u$ in $k$.

\medskip Assume that we are in case (a) : then by hypothesis $v(f(-1)) = v(f(1)) = 0$. For $v \in S_+ \cup S_-$,  we choose $\lambda_v$ such that

\medskip

$\underset {v \in S_+} \sum \partial [E_v,b_{\lambda_v},\alpha] = 0$ in $W(k) = W(k,\chi_+) \subset W_{\Gamma}(k),$ and

$\underset {v \in S_-} \sum \partial [E_v,b_{\lambda_v},\alpha] = 0$ in $W(k) = W(k,\chi_-) \subset W_{\Gamma}(k).$ 

\medskip \noindent
This
is possible by \cite {BT}, Proposition 6.6; indeed,  by \cite {BT}, Lemma 6.8 we have  $\underset {v \in S_=} \sum [\kappa_v : k] \equiv \ v(f(1)) \  ({\rm mod} \  2) $, 
$\underset {v \in S_-} \sum [\kappa_v : k] \equiv \ v(f(-1)) \  ({\rm mod} \  2) $, and $v(f(1)) = v(f(-1)) = 0$ by hypothesis;
therefore $\partial [E, b_{\lambda},\alpha] = 0$ in 
$W_{\Gamma}(k)$. 
Taking for $q^{\pm}$ the zero form if $n^{\pm} = 0$, and a unimodular form of determinant $u^{\pm}f(\pm 1)$ otherwise,  we get
$$\partial [E \oplus V, b_{\lambda} \oplus q, (\alpha,\epsilon)] = 0$$ in $W_{\Gamma}(k).$
This implies (i) and (ii), and completes the proof in case (a).

\bigskip Assume now that we are in case (b); then by hypothesis $v(f(-1)) = 0$ and $v(f(1)) = 1$. For $v \in  S_-$  we choose $\lambda_v$ such that

\medskip

$\underset {v \in S_-} \sum \partial [E_v,b_{\lambda_v},\alpha] = 0$ in $W(k) = W(k,\chi_-) \subset W_{\Gamma}(k).$ This
is possible by \cite {BT}, Proposition 6.6; indeed,  by \cite {BT}, Lemma 6.8 (ii) we have  $$\underset {v \in S_-} \sum [\kappa_v : k] \equiv \ v(f(-1)) \  ({\rm mod} \  2),$$ and $v(f(-1)) = 0$ by hypothesis. 

\medskip
We now come to the places in $S_+$. Recall that by \cite {BT}, Lemma 6.8 (i) we have $\underset {v \in S_+} \sum [\kappa_v : k] \equiv \ v(f(1)) \  ({\rm mod} \  2) $. Since $v(f(1)) = 1$ by hypothesis, this implies that $\underset {v \in S_+} \sum [\kappa_v : k] \equiv \ 1 \  ({\rm mod} \  2) $. 
Therefore there exists $w \in S_+$ such that $[\kappa_w:k]$ is odd. By \cite {BT}, Proposition 6.6, we can choose $\lambda_w$ such that
$\partial [E_w,b_{\lambda_w},\alpha]$ is either one of the two classes of $\gamma \in W(k) = W(k,\chi_+) \subset W_{\Gamma}(k)$ 
with ${\rm dim}(\gamma) = 1$. Let us choose the class of determinant $- \overline {u_+}$, and set $\partial [E_w,b_{\lambda_w},\alpha] = \delta$. 

\medskip Since $v(f(1)) = 1$, by hypothesis we have $n^+ \geqslant 1$. Let $(V^+,q^+)$ be a non-degenerate 
quadratic form over $K$ such that ${\rm det}(q^+) = u_+ f(1)$, and that 

\medskip
\centerline {$\partial[V^+,q^+,id] = - \delta$ in $W(k) = W(k,\chi_+) \subset W_{\Gamma}(k).$}

\medskip
Let $S_+' = S_+ - \{w \}$; we have $\underset {v \in S_+'} \sum [\kappa_v : k] \equiv \ 0 \  ({\rm mod} \  2) $, hence by \cite {BT}, Proposition 6.6,
for all $v \in S_+'$ there exists $\lambda_v \in E_{0,v}^{\times}$ such that 
$\underset {v \in S_+'} \sum \partial [E_v,b_{\lambda_v},\alpha] = 0$ in $W(k) = W(k,\chi_+) \subset W_{\Gamma}(k).$ 
We have

$$\partial [E \oplus V^+, b_{\lambda} \oplus q^+, (\alpha,id)] = 0$$ in $W_{\Gamma}(k)$.
Taking for $(V^-,q^-)$ a  quadratic form over $O$ of determinant $u_- f(-1)$ and setting $q = q^+ \oplus q^-$, we get
$$\partial [E \oplus V, b \oplus q, (\alpha,\epsilon)] = 0$$ in $W_{\Gamma}(k).$
This completes the proof in case (b). The proof is the same in case (c), exchanging the roles of $S_+$ and $S_-$. 

\bigskip
Assume that we are in case (d), that is, $v(f(1)) = v(f(-1)) = 1$. By \cite {BT}, Lemma 6.8 (i) and (ii), we have 

\medskip
$\underset {v \in S_+} \sum [\kappa_v : k] \equiv \ v(f(1)) \  ({\rm mod} \  2) $, and 
$\underset {v \in S_-} \sum [\kappa_v : k] \equiv \ v(f(-1)) \  ({\rm mod} \  2) $.  

\medskip Therefore
$\underset {v \in S_+} \sum [\kappa_v : k] \equiv \ \underset {v \in S_-} \sum [\kappa_v : k] \equiv \ 1 \  ({\rm mod} \  2) $. Hence
there exist $w_{\pm} \in S_{\pm}$ such that $[\kappa_{w_+}: k]$  and $[\kappa_{w_-}: k]$ are odd. By \cite {BT}, Proposition 6.6, we
can choose $\lambda_{w_{\pm}}$ such that $\partial [E_{w_{\pm}},b_{\lambda_{w_{\pm}}},\alpha]$ is either one of the two classes of $\gamma \in W(k) = W(k,\chi_{\pm}) \subset W_{\Gamma}(k)$ 
with ${\rm dim}(\gamma) = 1$. Let us choose  $\lambda_{w_{\pm}}$ such that  $\partial [E_{w_{\pm}},b_{\lambda_{w_{\pm}}},\alpha]$ is
represented by a form of dimension 1 and determinant $\overline u_{\pm}$, and set $$\delta_{\pm} = \partial [E_{w_{\pm}},b_{\lambda_{w_{\pm}}},\alpha].$$

\medskip By hypothesis, we have $n^+ \geqslant 1$ and $n^- \geqslant 1$. Let $(V^{\pm},q^{\pm})$ be non-degenerate 
quadratic forms over $K$ such that ${\rm det}(q^{\pm}) = u_{\pm}f(\pm 1)$ and that

\medskip
\centerline {$\partial[V^{\pm},q^{\pm},\epsilon^{\pm}] = - \delta_{\pm}$ in $W(k) = W(k,\chi_{\pm}) \subset W_{\Gamma}(k)$.}

\medskip
\medskip
Let $S_+' = S_+ - \{w_+ \}$; we have $\underset {v \in S_+'} \sum [\kappa_v : k] \equiv \ 0 \  ({\rm mod} \  2) $, hence by \cite {BT}, Proposition 6.6,
for all $v \in S_+'$ there exists $\lambda_v \in E_{0,v}^{\times}$ such that 
$\underset {v \in S_+'} \sum \partial [E_v,b_{\lambda_v},\alpha] = 0$ in $W(k) = W(k,\chi_+) \subset W_{\Gamma}(k).$ 
Similarly, set $S_-' = S_+ - \{w_- \}$; we have $\underset {v \in S_-'} \sum [\kappa_v : k] \equiv \ 0 \  ({\rm mod} \  2) $, hence by \cite {BT}, Proposition 6.6,
for all $v \in S_-$ there exists $\lambda_v \in E_{0,v}^{\times}$ such that 
$\underset {v \in S_-'} \sum \partial [E_v,b_{\lambda_v},\alpha] = 0$ in $W(k) = W(k,\chi_+) \subset W_{\Gamma}(k).$ 
Set $q = q^+ \oplus q^-$,
and note that $$\partial [E \oplus V, b_{\lambda} \oplus q, (\alpha,\epsilon)] = 0$$ in $W_{\Gamma}(k)$.
This completes the proof when the characteristic of $k$ is $ \not = 2$.

\medskip Assume now that the characteristic of $k$ is $2$. The algebra $E_v$ is of one of
the following types

\medskip
(sp) $E_v = E_{0,v} \times E_{0,v}$;

\medskip
(un) $E_v$ is an unramified extension of $E_{0,v}$;

\medskip
(r) $E_v$ is a ramified extension of $E_{0,v}$.

\medskip
This gives a partition $S = S_{sp} \cup S_{un} \cup S_r$.

\bigskip
If $\lambda = (\lambda_v)_{v \in S}$ is an element of $ E_0^{\times} =  \underset {v \in S} \prod E_{0,v}^{\times}$, 
note that by Lemma \ref{determinant} we have $${\rm disc}(b_{\lambda}) = (-1)^df(1)f(-1)$$ in $K^{\times}/K^{\times 2}$.

\medskip
We choose $\lambda = (\lambda_v)_{v \in S}$ in $E_0^{\times}$ such that 
for every $v \in S_{un}$, we have $\partial [E_v,b_{\lambda_v},\alpha] = 0$ in $W_{\Gamma}(k$); this is possible by \cite {BT}, Proposition 6.4. 

\medskip
Assume first that we are in case (a) or (d), and note that $v(f(1)) + v(f(-1)) = 0$. By \cite {BT}, Lemma 6.8 and Proposition 6.7, we can choose 
$\lambda_v$ such that 

\medskip

\centerline {$\underset {v \in S_r} \sum \partial [E_v,b_{\lambda_v},\alpha] = 0$ in $W(k) = W(k,1) \subset W_{\Gamma}(k).$}

\medskip Therefore
$\partial [E,b_{\lambda}, \alpha] = 0$ in $W_{\Gamma}(k)$.

\medskip
If we are in case (a) or case (d) and $K \not = {\bf Q}_2$, take for $q$ the unit quadratic form.  We have
$$\partial [E \oplus V, b_{\lambda} \oplus q, \alpha \oplus \epsilon] = 0$$ in $W_{\Gamma}(k)$; this concludes the proof in cases (a) and (d) when $K \not = {\bf Q}_2$.

\medskip
Suppose that 
$K = {\bf Q}_2$ and that we are in case (a). Suppose first that $n^+ = n^- = 0$. We already know that
$\partial [E,b_{\lambda}, \alpha] = 0$ in $W_{\Gamma}({\bf F}_2)$, hence (i) holds. Since $n^+ = n^- = 0$, by hypothesis $(-1)^df(1)f(-1) = 1$ in ${\bf Q}_2^{\times}/{\bf Q_2}^{\times 2}$, therefore ${\rm disc}(b_{\lambda}) = 1$ in ${\bf Q}_2^{\times}/{\bf Q_2}^{\times 2}$. 
Let us choose $\lambda$ such that the quadratic form
$(E,b_{\lambda})$ contains an even, unimodular 
${\bf Z}_2$-lattice. If $S_r = \varnothing$, this is automatic; indeed, in that case the trace map $E \to E_0$ is surjective, and hence
every ${\bf Z}_2$-lattice of the shape $(E,b_{\lambda})$  is even. If not, by  \cite {BT} Propositions  8.4 and 5.4 we can choose $\lambda$ having this additional property. This implies that (iii) holds as well. 

\medskip

We continue supposing that $K = {\bf Q}_2$ and that we are in case (a); assume now  that
$n^+ \not = 0$ and  $n^- = 0$. Let us choose $q^+$ such that ${\rm det}(q^*) = (-1)^n f(1) f(-1)$; since 
${\rm det}(b_{\lambda}) = f(1)f(-1)$, this implies that $${\rm det}(E \oplus V,b_{\lambda} \oplus q^+) = (-1)^n.$$ Moreover,
let us choose the Hasse-Witt invariant of $q^+$  in such a way that the quadratic form 
$(E \oplus V,b_{\lambda} \oplus q^+ \oplus q^-)$ contains an even, unimodular 
${\bf Z}_2$-lattice; this is possible by \cite {BT} Proposition 8.4. 
Therefore condition (iii) holds. Note that since $v({\rm det}(q^+)) = 0$, we have $\partial [V,q^+] = 0$ in $W({\bf F}_2$,
hence 
$\partial [E \oplus V, b_{\lambda} \oplus q, \alpha \oplus \epsilon] = 0$ in $W_{\Gamma}({\bf F}_2)$; therefore condition
(i) also holds. 

\medskip Assume now that $n^+ = 0$ and  $n^- \not = 0$. Let us choose $q^-$ such that ${\rm det}(q^-) = (-1)^n f(1) f(-1)$,
and note that since $v(f(1)) = v(f(-1)) = 0$, this implies that $v({\rm det}(q^-)) =  v(f(-1)) = 0$. As in the previous case,
we see that ${\rm det}(E \oplus V,b_{\lambda} \oplus q^-) = (-1)^n$, and we choose the Hasse-Witt invariant of
$q^-$ so that $(E \oplus V,b_{\lambda} \oplus q^-)$ contains an even, unimodular 
${\bf Z}_2$-lattice; this is possible by \cite {BT} Proposition 8.4. As in the previous case, we conclude that conditions (i)
and (iii) are satisfied. 

\medskip Suppose that $n^+ \not = 0$ and  $n^- \not = 0$.  Let us choose $q^+$ such that ${\rm det}(q^+) = u_+ f(1)$ and
$q^-$ such that ${\rm det}(q^-) = u_-  f(-1)$. Since $u_+ u_- = (-1)^n$ and 
${\rm det}(b_{\lambda}) = f(1)f(-1)$, this implies that $${\rm det}(E \oplus V,b_{\lambda} \oplus q^+ \oplus q^-) = (-1)^n.$$ Note
that since $v(u_-) = v(f(-1)) = 0$,  we have $v({\rm det}(q^-)) =  v(f(-1)) = 0$. As in the previous cases, we can
choose $q^+$ and $q^-$ such that $(E \oplus V,b_{\lambda} \oplus q^-)$ contains an even, unimodular 
${\bf Z}_2$-lattice, and that $\partial [E \oplus V, b_{\lambda} \oplus q, \alpha \oplus \epsilon] = 0$ in $W_{\Gamma}({\bf F}_2)$,
hence conditions (i)  and (iii) hold.

\medskip
Assume now that 
$K = {\bf Q}_2$ and that we are in case (d); note that the hypothesis implies that $n^+, n^- \geqslant 2$, and that both
$n^+$ and $n^-$ are even. With our previous choice of $\lambda$, we have $\partial [E,b_{\lambda}, \alpha] = 0$ in $W_{\Gamma}({\bf F}_2)$. Let us choose $q^+$ and $q^-$ such that ${\rm det}(q^{\pm}) = u_{\pm}f(\pm 1)$, and 
note that this implies that $v({\rm det}(q^-)) = v(f(-1))$,  and that  ${\rm det}(E \oplus V,b_{\lambda} \oplus q^+ \oplus q^-) = (-1)^n$. Moreover, choose the Hasse-Witt invariants of $q^+$ and $q^-$  so
that $(E \oplus V,b_{\lambda} \oplus q^+ \oplus q^-)$ contains an even, unimodular 
${\bf Z}_2$-lattice; this is possible by \cite {BT} Proposition 8.4. Therefore condition (iii) holds; moreover, we have
$\partial (V,q,\epsilon) = 0$ in $W_{\Gamma}({\bf F}_2)$ and
$\partial [E \oplus V, b_{\lambda} \oplus q, \alpha \oplus \epsilon] = 0$ in $W_{\Gamma}({\bf F}_2)$, 
hence condition (i) is also satisfied.  This concludes the proof in cases (a) and (d).

\medskip
Suppose that we are in case (b) or case (c), and note that in both cases, we have $v(f(1)) + v(f(-1)) = 1$. Hence Proposition 6.7 and Lemma 6.8 imply
that  $\underset {v \in S_r} \sum \partial [E_v,b_{\lambda_v},\alpha]$ is the unique non-trivial element of  $W(k) = W(k,1) \subset W_{\Gamma}(k).$ 
Suppose first that  $K \not = {\bf Q}_2$. We have either $n^+ \geqslant 1$ or $n^- \geqslant 1$; choose $q^{\pm}$ such that 
$\partial [V^{\pm},q^{\pm},\pm id]$ is also the unique non-trivial element of $W(k) = W(k,1) \subset W_{\Gamma}(k).$ 
We have $$\partial [E \oplus V, b_{\lambda} \oplus q, \alpha \oplus \epsilon] = 0$$ in $W_{\Gamma}(k)$. This settles cases (b) and (c)
when $K \not = {\bf Q}_2$.

\medskip
Assume now that $K = {\bf Q}_2$, and that we are in case (b), 
namely $v(f(1)) = 1$ and $v(f(-1)) = 0$; then $n^+ \ge 2$, and is even. If $n^- \not=  0$, then choose $q^-$ 
such that ${\rm det}(q^{-}) = u_{-}f(- 1)$, and note that this implies that $v({\rm det}(q^-)) = v(f(-1)) = 0$;
choose $q^+$ such that 
${\rm det}(q^{+}) = u_{+}f( 1)$. Since $v(f(1)) = 1$, this implies that  $\partial [V^+,q^+,id]$ is the unique non-trivial element of $W({\bf F}_2) = W_{\Gamma}({\bf F}_2,1) \subset W_{\Gamma}({\bf F}_2).$ Note that ${\rm det}(E \oplus V,b_{\lambda} \oplus q^+ \oplus q^-) = (-1)^n$ in  ${\bf Q}_2^{\times}/{\bf Q}^{\times 2}$. Moreover, choose the Hasse-Witt invariants of $q^+$ and $q^-$ such that the quadratic form
$(E \oplus V,b_{\lambda} \oplus q^+ \oplus q^-)$ contains an even, unimodular 
${\bf Z}_2$-lattice; this is possible by \cite {BT} Proposition 8.4. Hence condition  (iii) holds, and condition (i)
follows from the fact that $\partial [E,b_{\lambda}, \alpha]$ and $\partial [V^+,q^+,id]$ are both equal to the unique non-trivial element of $W({\bf F}_2) = W_{\Gamma}({\bf F}_2,1)$, which is a group of order $2$. Therefore 
$\partial [E \oplus V, b_{\lambda} \oplus q, \alpha \oplus \epsilon] = 0$ in $W_{\Gamma}({\bf F}_2)$, and
hence condition (i) is also satisfied.  

\medskip
Suppose now that $K = {\bf Q}_2$, and that we are in case  (c). Then $v(f(1)) = 0$ and $v(f(-1)) = 1$, hence  $n^- \ge 2$, and is even. If $n^+ \not=  0$, then choose $q^+$ such that ${\rm det}(q^{+}) = u_{+}f( 1)$. Choose 
$q^-$ 
such that ${\rm det}(q^{-}) = u_{-}f(- 1)$, and note that this implies that $v({\rm det}(q^-)) = v(f(-1)) = 1$, and that
$\partial [V^-,q^-,-id]$ is the unique non-trivial element of $W({\bf F}_2) = W_{\Gamma}({\bf F}_2,1) \subset W_{\Gamma}({\bf F}_2).$
We conclude as in case (b). 
This settles cases (b) and (c), and hence the
proof of the proposition is complete.

\bigskip We now show that the conditions of Theorem \ref{local odd} are sufficient. In the case where ${\rm char}(k) \not = 2$,
we obtain a more precise result (see part (ii) of the following result; for $K = {\bf Q}_2$ an analogous result
is given in Theorem \ref{local even}). We use the notation of \S  \ref{symmetric section}.

\begin{theo}\label{sufficient}  Let $F \in O[X]$ be a monic, symmetric polynomial.

\medskip If ${\rm char}(k) \not = 2$, assume that 
 $v(F(1)) = v(F(-1)) = 0$. 
 
\medskip If ${\rm char}(k) = 2$, assume that
 $v(F(1) F(-1)) = 0$. 
 
 \medskip
 {\rm (i)}
 Then there exists  a unimodular $O$-lattice having a
semi-simple isometry with characteristic polynomial $F$.

\medskip
{\rm (ii)} Assume in addition  that  ${\rm char}(k) \not = 2$ and let $u_+, u_- \in O^{\times}$. If $M^{\pm} \not = 0$, then there exists a unimodular  $O$-lattice 
having a semi-simple isometry with characteristic polynomial $F$ such that the associated $K[\Gamma]$-bilinear form
$(M^{\pm},q^{\pm})$ is such that  $${\rm det}(q^{\pm}) = u_{\pm} F_1(\pm 1)$$
in $K^{\times}/K^{\times 2}$.

\end{theo}

\noindent {\bf Proof.} Let us write $F = F_0 F_1 F_2$, where $F_i$ is the product of the irreducible factors of $F$ of type $i$.
The hyperbolic $O$-lattice of rank ${\rm deg}(F_2)$ has an isometry with characteristic polynomial $F_2$, therefore
it is enough to prove the theorem for $F = F_0 F_1$. 

 \medskip
From now on, we assume that $F = F_0 F_1$, in other words, all the irreducible factors of $F$ are symmetric, of type 0 or 1. Let $I_1$ be the set of irreducible factors
of type 1 of $F$. We have $F_1 = \underset{f \in I_1} \prod f^{n_f}$; note that $F_1(1)F_1(-1) \not = 0$. 
Let us write  $F(X) = F_1(X)(X-1)^{n^+}(X+1)^{n^-}$ for some integers $n^+, n^-$ such that $n^+,  n^- \geqslant 0$. 

\medskip
For all $f \in I_1$, set $E_f = K[X]/(f)$. Let $\sigma_f : E_f \to E_f$ be the involution induced by $X \mapsto X^{-1}$,
and let $(E_f)_0$ be the fixed field of $\sigma$ in $E_f$. Let $M_f$ be an extension of degree $n_f$ of $(E_f)_0$, linearly disjoint from $E_f$
over $(E_f)_0$. Set $\tilde E_f = E_f \otimes_{{(E_f)}_0} M_f$. Let $\alpha_f$ be a root of $f$ in $E_f$. The characteristic polynomial of the multiplication
by $\alpha_f$ on $\tilde E_f$ is $f^{n_f}$, and its minimal polynomial is $f$. Set $\tilde E = \underset {f \in I_1} \prod \tilde E_f$, and 
$\tilde M = \underset {f \in I_1} \prod \tilde M_f$. Let $\tilde \sigma : \tilde E \to \tilde E$ be the involution of $\tilde E$ induced by the involutions
$\sigma : E_f \to E_f$. Set $\tilde \alpha = (\alpha_f)_{f \in I_1}$, and let us denote by
$\tilde \alpha : \tilde E \to \tilde E$ the multiplication by 
$\tilde \alpha$. Note that $\tilde \alpha$ is semi-simple,
with characteristic polynomial $F_1$.

\medskip
 Let $V = V^+ \oplus V^-$ be a
$K$-vector space with ${\rm dim}(V^+) = n^+ $ and ${\rm dim}(V^-) = + n^-$. 
Applying Proposition \ref{technical} (i) with  $E_0 = \tilde M$, $E = \tilde E$, $\sigma = \tilde \sigma$,
$\alpha = \tilde \alpha$ and $f = F_1$, we see that 
there exists $\lambda \in \tilde M^{\times}$ and a non-degenerate quadratic form $q : V \times V \to K$ such that
$$\partial [\tilde E \oplus V, b_{\lambda} \oplus q, \alpha \oplus \epsilon] = 0$$ in $W_{\Gamma}(k)$,. By Theorem \ref{4.3} (iv) this implies that there exists a unimodular $O$-lattice having a semi-simple isometry with characteristic polynomial $F$,
proving part (i) of the theorem. Similarly, Proposition \ref{technical} (ii) implies part (ii) of the theorem.

 \bigskip
 To show that the conditions of Theorem \ref{local odd} are necessary, we start with some notation and a preliminary result.
 
 \medskip
 Extending the scalars to $K$, an even, unimodular lattice having a semi-simple isometry with characteristic polynomial $F$ gives
 rise to a $K[\Gamma]$-bilinear form on the semi-simple $K[\Gamma]$-module associated to $F$ (see \S \ref{symmetric section}), and
 this form has an orthogonal decomposition $M = M^0 \oplus M^1 \oplus M^2$, cf. \S \ref{isometries section}. The $K[\Gamma]$-form $M^0$
 has the further orthogonal decomposition $M^0 = (M^+,q^+) \oplus (M^-,q^-)$.

 \begin{notation}\label{N pm}  Let $\gamma$ be a generator of $\Gamma$. 
 Let $N_{\pm}$ be the simple $k[\Gamma]$-module such that ${\rm dim}_k(N_+) = 1$, and that $\gamma$ acts on $N_{\pm}$ by $\pm id$; note
 that $N_+ = N_-$ if  ${\rm char}(k) = 2$.
 
 \end{notation}
 
  \begin{lemma} \label{necessary lemma} Let $F \in O[X]$ be a monic, symmetric polynomial, and suppose that there exists a unimodular lattice having a semi-simple isometry
 with characteristic polynomial $F$. Let $M = M^0 \oplus M^1 \oplus M^2$ be the corresponding orthogonal decomposition of $K[\Gamma]$-bilinear forms. Let us write $F(X) = F_1(X)(X-1)^{n^+}(X+1)^{n^-}$ for some integers $n^+, n^-$ such that $n^+, n^- \geqslant 0$, and such that $F_1(1)F_1(-1) \not = 0$.  Then we have 
 
 \medskip
 {\rm (i)} Assume that ${\rm char}(k) \not = 2$. Then the component of $\partial[M^1]$ in $W_{\Gamma}(k,N_+) \simeq W(k)$ is represented by a 
 quadratic form of dimension $v(F_1(1))$ over $k$. 
 
  \medskip
 {\rm (i)} Assume that ${\rm char}(k) \not = 2$. Then the component of $\partial[M^1]$ in $W_{\Gamma}(k,N_-) \simeq W(k)$ is represented by a 
 quadratic form of dimension $v(F_1(-1))$ over $k$. 

 \medskip
 {\rm (i)} Assume that ${\rm char}(k) = 2$. Then the component of $\partial[M^1]$ in $W_{\Gamma}(k,N_+) = W_{\Gamma}(k,N_-)
  \simeq W(k)$ is represented by a 
 quadratic form of dimension $v(F_1(1))+v(F_1(-1))$ over $k$.

 \end{lemma}
 
 \noindent
 {\bf Proof.} Since $M$ is extended from a unimodular lattice, we have $\partial [M]  = 0$ (see Theorem \ref{4.3}). Let $M = M^0 \oplus M^1 \oplus M^2$ be the
 orthogonal decoposition of \S \ref{isometries section}. We have  $\partial [M^2] = 0$, hence $\partial ([M^0] + [M^1] )= 0$. 
 
 \medskip From now on, we assume that $M = M^0 \oplus M^1$; equivalently, all the irreducible factors of $F$ are of type 0 or 1. 
 Let us write  $F_1 = \underset {f \in I_1} \prod f^{n_f}$. 
 We have an orthogonal decomposition $$M^1 = \underset{f \in I} \oplus M_f,$$ where $M_f = [K[X]/(f)]^{n_f}$ (see \cite {M}, \S 3, or \cite {B 15},
 Propositions 3.3 and 3.4).  For all $f \in I$, set $E_f = K[X]/(f)$, and let $\sigma : E_f \to E_f$ be the $K$-linear involution induced by $X \mapsto X^{-1}$. 
 By a well-known transfer property (see for instance \cite {M}, Lemma 1.1 or \cite {B 15}, Proposition 3.6) the $K[\Gamma]$-bilinear form $M_f$
 is the trace of a non-degenerate hermitian form over $(E_f,\sigma)$, hence it is an orthogonal sum of forms of the type $b_{\lambda}$, see
 notation \ref{lambda}. 
 
 \medskip By \cite {BT}, Lemma 6.8 (i)
 and Proposition 6.6, 
the
 component of $\partial[M^1]$ in $W_{\Gamma}(k,M_+)$ is represented by a form of dimension $v(F_1(1))$, and this implies (i). Similarly, 
 applying  \cite {BT}, Lemma 6.8 (ii)
 and Proposition 6.6 implies (ii), and \cite {BT}, Lemma 6.8 (ii)
 and Proposition 6.7) implies (iii).

  \begin{prop} \label{necessary} Let $F \in O[X]$ be a monic, symmetric polynomial, and suppose that there exists a unimodular lattice having a semi-simple isometry
 with characteristic polynomial $F$. Then we have
 
  \medskip If ${\rm char}(k) \not = 2$, then  $v(F(1)) = v(F(-1)) = 0$.

  \medskip If ${\rm char}(k) = 2$, then $v(F(1) F(-1)) = 0$. 
   
 \end{prop}
 
 \noindent
 {\bf Proof.}  
 Suppose first that ${\rm char}(k) \not = 2$. 
If $n^+ > 0$ and $n^- > 0$, then $F(1) = F(-1) = 0$, so there is nothing to prove.  Assume that $n^+ = 0$. Then the component of $\partial [M^0]$ in $W_{\Gamma}(k,M_+)$ is trivial, and note that  $F(1) = F_1(1)$. By Lemma \ref{necessary lemma}  (i), the 
 component of $\partial[M^1]$ in $W_{\Gamma}(k,M_+)$ is represented by a form of dimension $v(F_1(1)) = v(F(1))$, hence  $v(F(1)) = 0$. 
 Similarly, $n^- = 0$ implies that $v(F(-1)) = 0$. This completes the proof of the proposition in the case where ${\rm char}(k) \not = 2$.
 
 \medskip
 Assume now that ${\rm char}(k) = 2$. If $n^+ > 0$ or $n^- > 0$, then $F(1)  F(-1) = 0$, so there is nothing to prove. Assume that $n^+ = n^- = 0$. The
 component of $\partial[M^1]$ in $W_{\Gamma}(k,M_+) = W_{\Gamma}(k,M_-)$ is represented by a form of dimension $v(F(1)) + v(F(-1))$ (cf. Lemma \ref{necessary lemma} (iii)). Since $n^+ = n^- = 0$ we have $M = M^1$, hence $\partial[M^1] = 0$,
 and we also have $F = F_1$; therefore $v(F(1)) + v(F(-1)) = 0$. This completes
 the proof of the proposition.
 
 \medskip
 \noindent {\bf Proof of Theorem \ref{local odd}.} The theorem follows from Theorem \ref{sufficient} and Proposition \ref{necessary}.

\section{Even, unimodular $\Gamma$-lattices over ${\bf Z}_2$}

We keep the notation of \S \ref{local}, with $K = {\bf Q}_2$ and $O = {\bf Z}_2$. Recall that if $a \in {\bf Z}_2$, we set
$v(a) = 0$ if $a = 0$ or if the $2$-adic valuation of $a$ is even, and $v(a) = 1$ if the $2$-adic valuation of $a$ is odd.

\medskip

If  $F \in {\bf Z}_2[X]$ is a monic, symmetric polynomial, we write $F = F_0 F_1 F_2$, where $F_i$ is the product of
the irreducible factors of type $i$ of $F$. Recall that $M = M^0 \oplus M^1 \oplus M^2$ is the semi-simple ${\bf Q}_2 [\Gamma]$-module
associated to $F$, and that $M_0 = M^+ \oplus M^-$.

\medskip
\begin{theo}\label {local even}  Let $F \in {\bf Z}_2[X]$ be a monic, symmetric polynomial of even degree such that  $F(0) = 1$, and set $2n = {\rm deg}(F)$. Let $u_+, u_- \in {\bf Z}_2^{\times}$ such that $u_+ u_- = (-1)^n$. Assume that the following conditions hold :

\medskip
{\rm (a)}  $v(F(1)) = v(F(-1)) = 0$.

\medskip
{\rm (b)} If $F(1)F(-1) \not = 0$, then  $(-1)^nF(1)F(-1) = 1$ in ${\bf Q}_2^{\times}/{\bf Q_2}^{\times 2}$.

\medskip
Then we have

\medskip
{\rm (i)} There  exists an even, unimodular ${\bf Z}_2$-lattice having a
semi-simple isometry with characteristic polynomial $F$.

\medskip
{\rm (ii)}  If $M^+ \not = 0$, $M^- \not = 0$ and $M^1 \not = 0$, then there exists an even,  unimodular  ${\bf Z}_2$-lattice 
having a semi-simple isometry with characteristic polynomial $F$ such that the associated ${\bf Q}_2[\Gamma]$-bilinear form
$(M^{\pm},q^{\pm})$ is such that  $${\rm det}(q^{\pm}) = u_{\pm} F_1(\pm 1)$$
in ${\bf Q}_2^{\times}/{\bf Q_2}^{\times 2}$.

\end{theo}

\noindent{\bf Proof.}  
Let $I$ be the set of irreducible factors of $F$ of type 1, and set
$F_1 = \underset{f \in I_1} \prod f^{n_f}$.
The hyperbolic ${\bf Z}_2$-lattice of rank ${\rm deg}(F_2)$ has a semi-simple isometry with characteristic polynomial $F_2$, therefore
it is enough to prove the theorem for $F = F_0 F_1$. 

 \medskip
From now on, we assume that all the irreducible factors of $F$ are symmetric, of type 0 or 1; we have $F = F_1(X-1)^{n^+}(X+1)^{n^-}$ for some integers $n^+  \geqslant 0$, $n^-  \geqslant 0$.
Note that since ${\rm deg}(F)$ is even and $F(0) = 1$ by hypothesis, $n^+$ and $ n^- $ are both even.

\medskip For all $f \in I_1$, set $E_f = {\bf Q}_2[X]/(f)$. Let $\sigma : E_f \to E_f$ be the involution induced by $X \mapsto X^{-1}$,
and let $(E_f)_0$ be the fixed field of $\sigma$ in $E_f$. Let $M_f$ be an extension of degree $n_f$ of $(E_f)_0$, linearly disjoint from $E_f$
over $(E_f)_0$. Set $\tilde E_f = E_f \otimes_{{(E_f)}_0} M_f$. Let $\alpha_f$ be a root of $f$ in $E_f$. The characteristic polynomial of the multiplication
by $\alpha_f$ on $\tilde E_f$ is $f^{n_f}$, and its minimal polynomial is $f$. Set $\tilde E = \underset {f \in I_1} \prod \tilde E_f$, and 
$\tilde M = \underset {f \in I_1} \prod \tilde M_f$. Let $\tilde \sigma : \tilde E \to \tilde E$ be the involution of $\tilde E$ induced by the involutions
$\sigma : E_f \to E_f$. Set $\tilde \alpha = (\alpha_f)_{f \in I_1}$, and let us denote by
$\tilde \alpha : \tilde E \to \tilde E$ the multiplication by 
$\tilde \alpha$. Note that $\tilde \alpha$ is semi-simple,
with characteristic polynomial $F_1$.

\medskip We apply Theorem 8.1 of \cite {BT} and Proposition \ref{technical} with  $E_0 = \tilde M$, $E = \tilde E$, $\sigma = \tilde \sigma$ and 
$\alpha = \tilde \alpha$. 

\medskip Let $V^{\pm}$ be a ${\bf Q}_2$-vector spaces of dimension $n^{\pm}$, and set $V = V^+ \oplus V^-$. Note that  
if $n^+ = n^- = 0$, then $F_1 = F$, hence the class of 
$(-1)^nF_1(1)F_1(-1) = 1$ in ${\bf Q}_2^{\times}/{\bf Q_2}^{\times 2}$ by hypothesis; therefore the hypotheses of  Proposition \ref{technical} 
are satisfied. 
Proposition \ref{technical} (i) and (iii) imply that there exist $\lambda \in \tilde M^{\times}$ and a non-degenerate quadratic form $q : V \times V \to {\bf Q}_2$
such that 
$$\partial [\tilde E \oplus V, b_{\lambda} \oplus q, \tilde {\alpha} \oplus \epsilon] = 0$$ in $W_{\Gamma}(k)$, that $(\tilde E \oplus V,b_{\lambda} \oplus q)$
contains an even, unimodular ${\bf Z}_2$-lattice, and that $v({\rm det}(q^-) = v(F_1(-1))$. By Theorem \ref{4.3} (iv), this
implies that there exists a unimodular lattice in $\tilde E \oplus V$ stable by $\tilde \alpha \oplus \epsilon$, hence a unimodular lattice having a semi-simple isometry with characteristic polynomial $F$; therefore conditions (i) and (ii)  of \cite {BT}, Theorem 8.1 hold. Moreover, since 
$v({\rm det}(q^-) = v(F_1(-1))$, Theorem 8.4 of \cite {BT} implies that condition (iii) of  \cite {BT}, Theorem 8.1 is also satisfied. 
This implies that there exists an even, unimodular ${\bf Z}_2$-lattice having a semi-simple isometry with characteristic polynomial $F$, and this completes the proof of (i). Part (ii) of the theorem  also follows from Proposition \ref{technical}, part (iii).

 \medskip
\begin{theo}\label {local even necessary}  Let $F \in {\bf Z}_2[X]$ be a monic, symmetric polynomial of even degree such that $F(0) = 1$, and set $2n = {\rm deg}(F)$. Assume that there
exists an even, unimodular ${\bf Z}_2$-lattice having a
semi-simple isometry with characteristic polynomial $F$. Then we have

\medskip
{\rm (a)}  $v(F(1)) = v(F(-1)) = 0$.

\medskip
{\rm (b)} If $F(1)F(-1) \not = 0$, then the class of $(-1)^nF(1)F(-1)$ in ${\bf Q}_2^{\times}/{\bf Q_2}^{\times 2}$ lies in $\{1,-3 \}$.

\end{theo} 

\noindent{\bf Proof.} 
 Let $L$ be an even, unimodular lattice  having a  semi-simple isometry with characteristic polynomial $F$.
 The lattice $L$ gives rise to a ${\bf Q}_2[ \Gamma]$-bilinear form $M$ on a bounded module. Let us consider
 the orthogonal decomposition of ${\bf Q}_2[ \Gamma]$-bilinear forms
 $$M = M^0 \oplus M^1 \oplus M^2$$ (cf. \S \ref{isometries section}). Since $L$ is unimodular we have $\partial [M]  = 0$; note that
 $\partial [M^2] = 0$, hence we have $\partial ([M^0] + [M^1] )= 0$. 
  
 \medskip From now on, we assume that $M = M^0 \oplus M^1$; equivalently, all the irreducible factors of $F$ are of type 0 or 1. 
 Let $I$ be the set of irreducible factors of $F$ of type 1, and set
$F_1 = \underset{f \in I_1} \prod f^{n_f}$.
 We have
$F = F_1(X-1)^{n^+}(X+1)^{n^-}$ for some integers $n^+  \geqslant 0$, $n^-  \geqslant 0$.

 \medskip
Further,  we have an orthogonal decomposition $M^1 = \underset{f \in I_1} \oplus M_f$, where $$M_f = [{\bf Q}_2[X]/(f)]^{n_f}$$ (see \cite {M}, \S 3, or \cite {B 15},
 Propositions 3.3 and 3.4). For all $f \in I_1$, set $E_f = {\bf Q}_2[X]/(f)$, and let $\sigma : E_f \to E_f$ be the ${\bf Q}_2$-linear involution induced by $X \mapsto X^{-1}$. 
 By a well-known transfer property (see for instance \cite {M}, Lemma 1.1 or \cite {B 15}, Proposition 3.6) the ${\bf Q}_2[\Gamma]$-bilinear form $M_f$
 is the trace of a non-degenerate hermitian form over $(E_f,\sigma)$, hence it is an orthogonal sum of forms of the type $b_{\lambda}$, see
 notation \ref{lambda}. 
 
 \medskip  The
 component of $\partial[M^1]$ in $W_{\Gamma}(k,N_{\pm})$ is represented by a form of dimension $v(F(1)) + v(F(-1))$ (mod 2) (cf. \cite {BT}, Lemma 6.8 (ii)
 and Proposition 6.7). 
 
 \medskip
 Suppose that $n^+ = n^- = 0$. Then $M = M^1$, hence we have  $\partial[M^1] =0$; by the above argument this implies that 
 $v(F(1)) + v(F(-1))$ (mod 2). By \cite {BT},  Proposition 8.6 and Theorem 8.5, we have $v(F(-1)) = 0$, hence (a) holds. Since $L$ is even and
 unimodular, the class of  $(-1)^n F(1) F(-1)$ in ${\bf Q}_2^{\times}/{\bf Q_2}^{\times 2}$ lies in 
 $\{1,-3 \}$; this shows that (b) holds as well.
 
 \medskip
 Let $M^0 = V^+ \oplus V^-$, and let $q^{\pm}$ the the quadratic form on $V^{\pm}$.
 
 \medskip
 Suppose that $n^+ \not = 0$, and $n^- = 0$. Then $F(1) = 0$, hence $v(F(1)) = 0$. Since $n^- = 0$, the quadratic form $q^-$ is the zero form,
 and $v({\rm det}(q^-) )=  0$. By \cite {BT}, Theorem 8.5 and Proposition 8.6, we have $v({\rm det}(q^-) ) = v(F(-1))$, hence $v(F(-1)) = 0$. This implies
 (a), and (b) is obvious since $F(1) = 0$.
 
 \medskip
 Assume now that $n^+ = 0$ and $n^- \not = 0$; then $F(-1) = 0$, hence (b) holds. By \cite {BT}, Theorem 8.5 and Proposition 8.6, we have $v({\rm det}(q^-)) = v(F(-1))$; since $F(-1) = 0$, this implies that $v({\rm det}(q^-)) = 0$. Therefore $\partial [M^0] = 0$. This implies that $\partial[M^1] = 0$, and hence 
 $v(F(1)) + v(F(-1)) = 0$. Since we already know that $v(F(-1)) = 0$, we obtain $v(F(1)) = 0$, and this implies (a).
 
 \medskip
 Finally, if $n^+ \not = 0$ and $n^- \not = 0$, then $F(1) = F(-1) = 0$, and hence both (a) and (b) hold. This concludes the proof of the theorem. 
  
  \section{Milnor signatures and Milnor indices}\label{Milnor}

The aim of this section is to recall from \cite {B 20} some notions of signatures and indices, inspired by Milnor's paper \cite {M 68}.

\medskip Let $F \in {\bf R}[X]$ be a monic, symmetric polynomial.
If  $(V,q)$ is a non-degenerate quadratic form over $\bf R$ and if $t : V \to V$ is a semi-simple isometry of $q$ with characteristic
polynomial $F$, we associate to each irreducible, symmetric factor $\mathcal P$ of $F$ an
index $\tau({\mathcal P})$ and a signature $\mu({\mathcal P})$ as follows. Let $V_{\mathcal P(t)}$ be the $\mathcal P(t)$-primary subspace of $V$, consisting
of all $v \in V$ with $\mathcal P(t)^N v = 0$ for $N$ large. The {\it Milnor index} $\tau(\mathcal P)$ is by definition the index
of the restriction of $q$ to the subspace $V_{\mathcal P(t)}$,
and we define the {\it Milnor signature} $\mu({\mathcal P})$ at $\mathcal P$ as the signature of the restriction of $q$ to $V_{\mathcal P(t)}$. 

\medskip
Let  ${\rm Irr}_{\bf R}(F)$ be the set of irreducible, symmetric factors 
of $F  \in {\bf R}[X]$; if $\mathcal P \in {\rm Irr}_{\bf R}(F)$, then either ${\rm deg}({\mathcal P}) = 2$, or ${\mathcal P}(X) = X \pm 1$.  If $(r,s)$ is the signature of $q$, we have
$$\underset {\mathcal P} \sum \tau({\mathcal P}) = r-s,$$
where the sum runs over $\mathcal P \in {\rm Irr}_{\bf R}(F)$.

\medskip

If $\mathcal P \in {\rm Irr}_{\bf R}(F)$, let  $n_{\mathcal P} > 0$ be the integer such that
$\mathcal P^{n_{\mathcal P}}$ is the power of $\mathcal P$ dividing $F$.

\medskip We denote by ${\rm Mil}(F)$ the set of maps $\tau : {\rm Irr}_{\bf R}(F) \to {\bf Z}$ such that the image of 
$\mathcal P \in {\rm Irr}_{\bf R}(f)$ belongs to the set $\{ -{\rm deg}({\mathcal P})n_{\mathcal P},\dots,{\rm deg}({\mathcal P})n_{\mathcal P} \}$. 
For all integers $r, s \geqslant 0$, let
${\rm Mil}_{r,s}(F)$ be the subset of ${\rm Mil}(F)$ such that $\underset {\mathcal P} \sum \tau({\mathcal P}) = r-s$,
where the sum runs over $\mathcal P \in {\rm Irr}_{\bf R}(F)$.

\begin{prop}\label{bijection} Sending a semi-simple element ${\rm SO}_{r,s}({\bf R})$ with characteristic polynomial $F$ to its Milnor index
induces a bijection between the conjugacy classes of semi-simple elements of  ${\rm SO}_{r,s}({\bf R})$ and ${\rm Mil}_{r,s}(F)$.

\end{prop}

\noindent {\bf Proof.} See  \cite {B 20}, \S 6. 
 
 \section{Local conditions for even, unimodular $\Gamma$-lattices}\label{local section}
 
 Let $F \in {\bf Z}[X]$ be a monic, symmetric polynomial, and let
 $r, s \ge 0$ be integers such that $r + s = {\rm deg}(F)$.  The aim of this section is to give local conditions for the existence of an even, unimodular lattice of signature $(r,s)$ having a semi-simple
isometry with characteristic polynomial $F$. More precisely, given a Milnor index $\tau \in {\rm Mil}_{r,s}(F)$, we give a necessary and
sufficient condition for an even, unimodular lattice having a semi-simple isometry with characteristic polynomial $F$ 
and Milnor index $\tau$ to exist everywhere locally.

\medskip
Let $m(F)$ be the number of roots  $z$ of $F$ with $|z| > 1$ (counted with multiplicity). 

\begin{prop}\label{reals} There exist an ${\bf R}$-vector space $V$ and a non-degenerate quadratic form $q$ of signature $(r,s)$ having
a semi-simple isometry with characteristic polynomial $F$ and Milnor index $\tau \in {\rm Mil}_{r,s}(F)$ if and only if 
$r \geqslant m(F)$, $s \geqslant m(F)$, and if moreover $F(1)F(-1) \not = 0$, then $m(F) \equiv r \equiv s  \ {\rm (mod \ 2)}$.

\end{prop}

\noindent
{\bf Proof.} This follows from \cite {B 15}, Proposition 8.1. Indeed, the necessity of the conditions follows immediately from 
 \cite {B 15}, Proposition 8.1. To prove the sufficiency, note that while the statement of  \cite {B 15}, Proposition 8.1 only claims
 the existence of a non-degenerate quadratic form $q$ of signature $(r,s)$ having
a semi-simple isometry with characteristic polynomial $F$, the proof shows the existence of such a form
having a semi-simple isometry with a given Milnor index $\tau \in {\rm Mil}_{r,s}(F)$.

\bigskip
 If $p$ is a prime number, we say
that a ${\bf Z}_p$-lattice $(L,q)$ is even if $q(x,y) \in 2{\bf Z}_p$; note that if $p \not = 2$, then every lattice is even, since $2$
is a unit in ${\bf Z}_p$. 

\medskip
Assume moreover that $F$ has even degree, and that $F(0) = 1$. Set $2n = {\rm deg}(F)$. 

\begin{theo}\label{integral local}  There exists an even,  unimodular ${\bf Z}_p$-lattice having a semi-simple isometry
with characteristic polynomial $F$  for all prime numbers $p$ if and only if  $|F(1)|$, $|F(-1)|$ and $(-1)^n F(1) F(-1)$ are all squares. 

\end{theo}

\noindent
{\bf Proof.} This follows from Theorems \ref{local odd},  \ref{local even} and \ref{local even necessary}. Indeed,  if $p$ is a prime number $ \not = 2$, the existence of a unimodular 
${\bf Z}_p$-lattice having a semi-simple isometry
with characteristic polynomial $F$ implies that either $F(1) = 0$, or $v_p(F(1))$ is even; similarly, either $F(-1) = 0$, or $v_p(F(-1))$ is even
(see Theorem \ref{local odd}). The existence of an even, unimodular ${\bf Z}_2$-lattice implies the same property for $p = 2$ by 
Theorem \ref{local even necessary}. This implies that $|F(1)|$ and $|F(-1)|$ are both squares, and therefore $|F(1)F(-1)|$ is a square. If $F(1)F(-1) = 0$,
we are done. If not,
Theorem \ref{local even necessary} implies that the class of  $(-1)^n F(1) F(-1)$ in ${\bf Q}_2^{\times}/{\bf Q_2}^{\times 2}$ lies in 
 $\{1,-3 \}$; since $|F(1)F(-1)|$ is a square, this implies that
 $(-1)^n F(1) F(-1)$ is a square. The converse
is an immediate consequence of Theorems \ref{local odd} and \ref{local even}.




\section{The local-global problem}\label{local-global problem section}

The aim of this section is to reformulate the local conditions of \S \ref{local section}, and to give a framework for the local-global problem
of the next sections. We also introduce some notation that will be used in the following sections, 

\medskip

Let $F \in {\bf Z}[X]$ be a monic, symmetric polynomial of even degree such that $F(0) = 1$; set $2n = {\rm deg}(F)$. 
Let $J$ be the set of irreducible factors of $F$, and let us write $F = \underset{f \in J} \prod f^{n_f}$. Let $I_1 \subset J$ be the subset of irreducible factors of type 1, and let $I_0 \subset J$ be the
set of irreducible factors of type 0. 

\medskip Let $M = M^0 \oplus M^1 \oplus M^2$ be the semi-simple ${\bf Q}[\Gamma]$-module associated to the polynomial $F$ (see \S
\ref{symmetric section}).

\bigskip
Let $r, s \ge 0$ be integers such that $r + s = {\rm deg}(F)$ and that  $r  \equiv s \ {\rm (mod \ 8)}$.  Let $(V,q) = (V_{r,s},q_{r,s})$ be
the diagonal quadratic form over ${\bf Q}$ with $r$ entries $1$ and $s$ entries $-1$. 

\begin{prop}\label{equivalent local} The following properties are equivalent

\medskip
{\rm (i)} For all prime numbers $p$ there exists an even, unimodular ${\bf Z}_p$-lattice having a semi-simple isometry with characteristic
polynomial $F$.

\medskip
{\rm (ii)} $|F(1)|$, $|F(-1)|$ and $(-1)^n F(1) F(-1)$ are all squares. 

\medskip
{\rm (iii)} For all prime numbers $p$, the quadratic form $(V,q) \otimes_{\bf Q} {\bf Q}_p$ has a semi-simple isometry with characteristic
polynomial $F$ that stabilizes an even, unimodular lattice.

\medskip
{\rm (iv)}  For all prime numbers $p$, the quadratic form $(V,q) \otimes_{\bf Q} {\bf Q}_p$ has an isometry with module $M \otimes _{\bf Q} {\bf Q}_p$,
giving rise to a class $[M \otimes _{\bf Q} {\bf Q}_p,q]$ in $W_{\Gamma}({\bf Q}_p)$ such that $\partial_p [M \otimes _{\bf Q} {\bf Q}_p,q] = 0$
in $W_{\Gamma}({\bf F}_p)$ and that $v_2({\rm det}(q_-)) \equiv v_2(F_1(-1)) \ {\rm (mod \ 2)}$. 

\end{prop}

\noindent
{\bf Proof.} The equivalence of {\rm (i)} and {\rm (ii)} is Proposition \ref{integral local}, and it is clear that {\rm (iii)} implies {\rm (i)}. Let us show
that {\rm (i)} implies {\rm (iii)}. Set $u = (-1)^s$. Since $r + s = 2n$ and $r  \equiv s \ {\rm (mod \ 8)}$, we have $n  \equiv s \ {\rm (mod \ 8)}$, hence $u = (-1)^n$. 
 By {\rm (i)}, there exists an
even, unimodular ${\bf Z}_p$-lattice having a semi-simple isometry with characteristic polynomial $F$. If $F(1)F(-1) \not = 0$, then
the class of the determinant of this lattice in in ${\bf Q}_p^{\times}/{\bf Q_p}^{\times 2}$ is equal to $F(1)F(-1)$, and  $F(1)F(-1) = u$ by
{\rm (ii)}; if $F(1)F(-1)  = 0$,
then by Theorem \ref{sufficient} (ii) and Theorem \ref{local even} (ii) we can assume that the determinant of this lattice  in ${\bf Q}_p^{\times}/{\bf Q_p}^{\times 2}$ is equal to $u$. Therefore the
lattice is isomorphic to the diagonal ${\bf Z}_p$-lattice $\langle 1,\dots, u \rangle$ of determinant $u$ if $p \not = 2$ (cf.  \cite{OM}, 92:1), and to the orthogonal sum of $n$ hyperbolic planes if $p = 2$ (see for instance \cite {BT}, Proposition 8.3). 
Let $q^p$ be the
quadratic form over ${\bf Q}_p$ obtained from this lattice by extension of scalars; then the Hasse-Witt invariant of $q^p$ is trivial
if $p \not = 2$, and is equal to the Hasse-Witt invariant of the orthogonal sum of $n$ hyperbolic planes if $p = 2$; its 
determinant is equal to $u = (-1)^n$ in ${\bf Q}_p^{\times}/{\bf Q_p}^{\times 2}$. This implies that   $q^p$ and
$(V,q) \otimes {\bf Q}_p$ are isomorphic as quadratic forms over ${\bf Q}_p$. Since $q^p$ has a semi-simple isometry with characteristic
polynomial $F$ that stabilizes a unimodular lattice, property (iii) holds.
Finally, the equivalence of (iii) and (iv) follows from Theorem \ref{4.3} (iv) and from \cite{BT}, Theorems 8.1 and 8.5. 

\medskip
\noindent
{\bf Terminology.} We say that the {\it local conditions for $F$} hold at the finite places if the equivalent conditions of Proposition \ref{equivalent local}
are satisfied.

\bigskip
Recall that $m(F)$ is the number of roots  $z$ of $F$ with $|z| > 1$ (counted with multiplicity). 

\begin{prop}\label{equivalent real} Let $\tau \in {\rm Mil}_{r,s}(F)$ be a Milnor index. The following properties
are equivalent :

\medskip
{\rm (i)} The quadratic form $(V,q) \otimes_{\bf Q} {\bf R}$ has a semi-simple isometry with characteristic polynomial $F$ and
Milnor index $\tau$. 

\medskip
{\rm (ii)} $r \geqslant m(F)$, $s \geqslant m(F)$, and if moreover $F(1)F(-1) \not = 0$, then $m(F) \equiv r \equiv s  \ {\rm (mod \ 2)}$.

\medskip
{\rm (iii)} The quadratic form $(V,q) \otimes_{\bf Q} {\bf R}$ has an isometry with module $M \otimes_{\bf Q} {\bf R}$ and 
Milnor index $\tau$.

\end{prop}

\noindent {\bf Proof.} The equivalence of (i) and (ii) follows from Proposition \ref{reals}, and (iii) is a reformulation of (i). 

\medskip
\noindent
{\bf Terminology.} We say that the {\it local conditions for $F$ and $\tau$} hold at the infinite place if the equivalent conditions of Proposition \ref{equivalent real}
are satisfied.

\medskip We consider the following conditions

\medskip

(C 1) {\it $|F(1)|$, $|F(-1)|$ and $(-1)^n F(1) F(-1)$ are all squares.}

\medskip

(C 2) {\it $r \geqslant m(F)$, $s \geqslant m(F)$, and if moreover $F(1)F(-1) \not = 0$, then $m(F) \equiv r \equiv s  \ {\rm (mod \ 2)}$.}

\medskip
Note that the local conditions for $F$ at the finite places hold if and only if condition (C 1) is satisfied, and that the local conditions
for $F$ and $\tau$ hold if and only if condition (C 2) is satisfied.

\medskip
\noindent
{\bf Terminology.} Let $M$ and $q$ be as above, and let $p$ be a prime number. A $\Gamma$-quadratic form $(M \otimes _{\bf Q} {\bf Q}_p,q)$ such that $\partial_p[(M \otimes _{\bf Q} {\bf Q}_p,q] = 0$ in $W_{\Gamma}({\bf F}_p)$ 
and that $v_2({\rm det}(q_-)) \equiv v_2(F_1(-1)) \ {\rm (mod \ 2)}$ if $p = 2$ 
is called a {\it local solution} for $F$ at the prime number $p$. 

\medskip 

\section{${\bf Q}[\Gamma]$-forms, signatures and determinants}\label{Q section}

Let $F \in {\bf Z}[X]$ be a monic, symmetric polynomial, and let us write $F = F_0 F_1 F_2$, where $F_i$ is the product
of the irreducible factors of type $i$ of $F$. Let $r, s \ge 0$ be integers such that $r + s = {\rm deg}(F)$ and that  $r  \equiv s \ {\rm (mod \ 8)}$, and let $\tau \in {\rm Mil}_{r,s}(F)$ be a Milnor index. Let $(L,q)$ be an even, unimodular
lattice having a semi-simple isometry with characteristic polynomial $F$ and Milnor index $\tau$, and let $(M,q)$ be
the corresponding ${\bf Q}[\Gamma]$-form, and let $$M = M^0 \oplus M^1 \oplus M^2$$ and $$M^0 = M^+ \oplus M^-$$
be the associated orthogonal decompositions (cf. \S \ref{isometries section}). Note that the Milnor index $\tau$ and the degrees
of the polynomials determine the signatures of the factors. We have
${\rm sign}(M) = (r,s)$. Set ${\rm sign}(M^1) = (r_1,s_1)$, and ${\rm sign}(M^2) = (r_2,s_2)$; note that $r_2 = s_2 = {\rm deg}(F_2)/2$, since $M^2$ is hyperbolic, and set $${\rm sign}(M^+) = (r^+,s^+), \ \  {\rm sign}(M^-) = (r_-,s_-).$$ 

\smallskip
\noindent
We have ${\rm det}(M) = (-1)^s$, ${\rm det}(M^1) = F_1(1)F_1(-1) = (-1)^{s_1} |F_1(1)F_1(-1)|$, and ${\rm det}(M^2) = (-1)^{s_2}$. 

\begin{prop}\label{det} We have $${\rm det}(M^+) = (-1)^{s_+}|F_1(1)|, \ \  {\rm det}(M^-) = (-1)^{s_-}|F_1(-1)|$$ in ${\bf Q}^{\times}/{\bf Q}^{\times 2}$. 

\end{prop}

\noindent
{\bf Proof.} The sign of ${\rm det}(M^{\pm})$ is $(-1)^{s_{\pm}}$. Let $p$ be a prime, $p \not = 2$. 
By Lemma \ref{necessary lemma}, 
the component of $\partial_p[M^1]$ in $W_{\Gamma}({\bf F}_p,N_{\pm}) \simeq W({\bf F}_p)$ is represented by a 
 quadratic form of dimension $v(F_1(\pm 1))$ over $\bf F_p$. If $v(F_1(\pm 1)) = 0$, then this component
 of $\partial_p[M^1]$  is trivial, hence $\partial_p[M^{\pm}]$ is also trivial. This implies that $v({\rm det}(M^{\pm}) = 0$.
 Assume now that  $v(F_1(\pm 1)) = 1$. Then the component of $\partial_p[M^1]$ in $W_{\Gamma}({\bf F}_p,N_{\pm}) \simeq W({\bf F}_p)$ is represented by a 
 quadratic form of dimension 1. Since $\partial_p[M] = 0$, this implies that $\partial_p[M^{\pm}]$ is represented by a form
 of dimension 1, and therefore $v({\rm det}(M^{\pm}) = 1$. Hence in this case too, we have 
 $v({\rm det}(M^{\pm}))  = v(F_1(\pm 1))$. 

\medskip

Assume that $p = 2$. The component of $\partial_2[M^1]$ in $W_{\Gamma}({\bf F}_2,N_{\pm}) \simeq W({\bf F}_2)$ is represented by a 
 quadratic form of dimension $v(F_1(1)) + v(F_1(-1))$ over ${\bf F}_2$ (see Lemma \ref{necessary lemma}). If $M^+ = 0$
 and $M^- = 0$, there is nothing to prove. Assume that $M^+ \not = 0$, and $M^- = 0$. Then $F_1(-1) = F(-1)$, and
 by Theorem \ref{local even necessary} (a), we have $v(F_1(-1)) = 0$. Hence $\partial_2(M^1)$ is represented by a form of dimension 
 $v(F_1(1))$. If $v(F_1(1)) = 0$, then $\partial_2(M^1)=0$, and hence $\partial_2(M^+) = 0$; therefore $v({\rm det}(M^+)) = 0$. 
 If $v(F_1(1)) = 1$, then $\partial_2(M^1)$ is represented by a form of dimension 1 over ${\bf F}_2$, hence 
 $\partial_2(M^+)$ is also represented by a form of dimension 1 over ${\bf F}_2$. This implies that $v({\rm det}(M^+)) = 1$. 
 Therefore $v({\rm det}(M^+)) = v(F_1(1))$. The same argument shows that if $M^+ = 0$ and $M^- \not = 0$, then 
 $v({\rm det}(M^-)) = v(F_1(-1))$.  Suppose now that $M^+ \not = 0$ and $M^- \not = 0$. By \cite {BT}, Theorem 8.5 and
 Proposition 8.6, we have $v({\rm det}(M^-)) = v(F_1(-1))$. If $v(F_1(1)) = v(F_1(-1))$, then  $\partial_2(M^1) = 0$. 
 Therefore $\partial_2(M^+ \oplus M^-) = 0$. 
 Since  $v({\rm det}(M^-)) =  v(F_1(-1))$, this implies that $v({\rm det}(M^+)) = v(F_1(1))$. If 
 $v(F_1(1)) + v(F_1(-1))= 1$, then $\partial_2(M^1) \not = 0$, and hence $\partial_2(M^+ \oplus M^-) \not = 0$. 
 Therefore  $v({\rm det}(M^+)) + v({\rm det}(M^-))  = 1$. Since $v({\rm det}(M^-)) = v(F_1(-1))$, we have
 $v({\rm det}(M^+)) = v(F_1(1))$. This completes the proof of the proposition. 

\bigskip

\section{Local decomposition}\label{decomposition section}

Let $F \in {\bf Z}[X]$ be a monic, symmetric polynomial of even degree with $F(0) = 1$; set $2n = {\rm deg}(F)$. 
Let $r, s \ge 0$ be integers such that $r + s = {\rm deg}(F)$ and that  $r  \equiv s \ {\rm (mod \ 8)}$, let $\tau \in {\rm Mil}_{r,s}(F)$ be a Milnor index. If  the local conditions (C 1) and (C 2) hold, then we obtain
a local solution everywhere (see  \S \ref{local-global problem section}). The aim of this section is to define
local decompositions that will be useful in the following sections. 

\medskip We start by introducing some notation.
Let
$M = M^0 \oplus M^1 \oplus M^2$ be the semi-simple ${\bf Q}[\Gamma]$-module associated to $F$ as in \S \ref{symmetric section}, with 

\medskip
\centerline {$M^1 = \underset{i \in I} \oplus M_f$ and $M^0 = M^+ \oplus M^-$.} 

\medskip If $f \in I_1$, set $E_f = {\bf Q}[X]/(f)$ and let $\sigma_f : E_f \to E_f$ be the involution induced by $X \mapsto X^{-1}$.
Let $(E_f)_0$ be the fixed field of $\sigma_f$, and let $d_f \in (E_f)_0$ be such that $E_f = (E_f)_0 (\sqrt d_f)$. 
Note that $M_f$ is an $E_f$-vector space of dimension $n_f$. 
Let $$Q_f : M_f \times M_f \to {\bf Q}$$ be the orthogonal sum of $n_f$ copies of the quadratic form $E_f \times E_f \to {\bf Q}$ defined by $(x,y) \mapsto {\rm Tr}_{E_f/{\bf Q}}(x \sigma_f(y))$.

\medskip The Milnor index $\tau \in {\rm Mil}_{r,s}(F)$ determines the signatures of $M^+$ and $M^-$, as follows. 
Recall that ${\rm dim}(M^+) = n_+$ and ${\rm dim}(M^+) = n_+$.

\medskip Let $s_+$ and $s_-$ be
as in \S \ref{Q section}, and set $D_{\pm} = (-1)^{s_{\pm}}F_1(\pm 1)$; let $Q{\pm}$ be the diagonal quadratic form 
of dimension $n_{\pm}$ over ${\bf Q}$
defined by $Q{\pm} = \langle D_{\pm}, 1,\dots,1 \rangle$.
Let $Q$ be the orthogonal sum $$Q = \underset{f \in I_1} \oplus Q_f \oplus Q_+ \oplus Q_-.$$

\medskip

 $${\rm sign}(M^+) = (r^+,s^+), \ \  {\rm sign}(M^-) = (r_-,s_-).$$ 
 
\medskip

\medskip Recall form \S \ref{local-global problem section} that we denote by $q = q_{r,s}$ the diagonal
quadratic form over ${\bf Q}$ having $r$ diagonal entries $1$ and $s$ diagonal entries $-1$.

\medskip Assume that conditions (C 1) and (C 2) hold. If $p$ is a prime number, then $(M,q) \otimes_{\bf Q} {\bf Q}_p$ has
a structure of a ${\bf Q}_p[\Gamma]$-quadratic form (see \S \ref{local-global problem section}), and we have the
orthogonal decomposition (cf. \S \ref{isometries section}).

$$(M,q) \otimes_{\bf Q} {\bf Q}_p = \underset {f \in I_1} \oplus (M^p_f,q_f^p) \oplus (M^p_+,q_+^p) \oplus (M^p_-,q_-^p) \oplus (M^p_2,q_2^p),$$

\noindent
where $M^p_f = M_f \otimes_{\bf Q} {\bf Q}_p$, $M_+^p = M^+ \otimes_{\bf Q} {\bf Q}_p$, $M_-^p = M^- \otimes_{\bf Q} {\bf Q}_p$, and $M_2^p = M^2 \otimes_{\bf Q} {\bf Q}_p$. The ${\bf Q}_p[\Gamma]$-quadratic form $(M^p_2,q_2^p)$ is hyperbolic.

\medskip
 For $f \in I_1$, set
$E_f^p = E_f \otimes_{\bf Q} {\bf Q}_p$ and $(E_f)_0^p = (E_f)_0 \otimes_{\bf Q} {\bf Q}_p$. There exists a unique non-degenerate hermitian form $(M_f^p,h^p_f)$ over
$(E_f^p,\sigma_f)$ such that $$q_f^p(x,y) = {\rm Tr}_{E_f^p/{\bf Q}_p}(h^p_f(x,y)),$$
see for instance \cite {M}, Lemma 1.1 or \cite {B 15}, Proposition 3.6.

\medskip
Set 
$\lambda^p_f = {\rm det}(h^p_f) \in (E^p_f)_0^{\times}/{\rm N}_{E^p_f/(E^p_f)_0}$. Note that the hermitian
form $h^p_f$ is isomorphic to the $n_f$-dimensional diagonal hermitian form $\langle \lambda^p_f,1,\dots,1 \rangle$
over $E_f^p$. Hence $q_f^p$ is determined by $\lambda^p_f$.

\begin{notation}\label{lambda bis} With the notation above, set 
$$\partial_p(\lambda^p_f) = \partial_p[q_f^p] \in W_{\Gamma}({\bf F}_p).$$

\end{notation}

\smallskip

\begin{prop}\label{dim det w f} We have ${\rm dim}(q_f^p) = {\rm deg}(f)n_f$, ${\rm det}(q_f^p) = [f(1)f(-1)]^{n_f}$, and
the Hasse-Witt invariant of $q_f^p$ satisfies
$$w_2(q_f^p) + w_2(Q_f) = {\rm cor}_{(E_f)_0^p/{\bf Q}_p}({\rm det}(h^p_f),d_f)$$ in ${\rm Br}_2({\bf Q}_p)$. 

\end{prop}

\noindent
{\bf Proof.} The assertion concerning the dimension is clear, the one on the determinant follows from 
Lemma \ref{determinant}, and the property of the Hasse-Witt invariants from \cite {B 20}, Proposition 12.8.

\begin{prop}\label{dim det w pm} 
{\rm (a)} ${\rm dim}(q_{\pm}^p) = n_{\pm}$ and $${\rm det}(q_{+}^p){\rm det}(q_{+}^p) = (-1)^{s_+ + s_-}|F_1(1)F_1(-1)|.$$

\medskip
\noindent
{\rm (b)} If $n_+ \not = 0$ and $n_- \not = 0$, then we can choose $q_+^p$ and $q_-^p$ such that 
${\rm det}(q_{\pm}^p) = (-1)^{s_{\pm}}|F_1(\pm 1)|$. 

\medskip 
\noindent
{\rm (c)} If $n_{\pm} \not = 0$, then the 
Hasse-Witt invariant of $q_{\pm}^p$ can take either of the two possible values of
$\{0,1\} = {\rm Br}_2({\bf Q}_p)$. 

\medskip
\noindent
{\rm (d)} If $p \not = 2$, then $\partial_p[q_{\pm}^p]$ can be either of the
two possible classes of dimension $v_p({\rm det}(q_{\pm}^p))$ of $W_{\Gamma}({\bf F}_p,N^{\pm}) \simeq W({\bf F}_p)$.

\end{prop}

\noindent 
{\bf Proof.} (a) is clear, (b) follows from Theorem \ref{sufficient} (ii) and Theorem \ref{local even} (ii); (c) and (d) are 
straightforward to check.

\medskip We also need the following 

\begin{lemma} Let $p$ be a prime number, $p \not = 2$, and let $b_1$ and $b_2$ be two quadratic forms over ${\bf Q}_p$
with ${\rm dim}(b_1) = {\rm dim}(b_2)$ and ${\rm det}(b_1) = {\rm det}(b_2)$. Then we have
$$w_2(b_1) = w_2(b_2) \ \ {\rm in} \ \ {\rm Br}_2({\bf Q}_p) \ \ \iff \ \ \partial_p[b_1] = \partial_p[b_2] \ \ {\rm in} \ \ W({\bf F}_p).$$

\end{lemma}

\noindent
{\bf Proof.} The proof is straightforward.

\bigskip
Similarly, we have
$$(M,q) \otimes_{\bf Q} {\bf R} = \underset {f \in I_1} \oplus (M^{\infty}_f,q_f^{\infty}) \oplus (M^{\infty}_+,q_+^{\infty}) \oplus (M^{\infty}_-,q_-^{\infty}) \oplus (M^{\infty}_2,q_2^{\infty}),$$

\noindent
where $M^{\infty}_f = M_f \otimes_{\bf Q} {\bf R}$, $M_+^{\infty} = M^+ \otimes_{\bf Q} {\bf R}$, $M_-^{\infty} = M^- \otimes_{\bf Q} {\bf Q}_p$, and $M_2^{\infty} = M^2 \otimes_{\bf Q} {\bf Q}_p$. The ${\bf R}[\Gamma]$-quadratic form $(M^{\infty}_2,q_2^p)$ is hyperbolic.

\medskip
The  ${\bf R}[\Gamma]$-quadratic forms $(M^{\infty}_f,q_f^{\infty})$ and $(M^{\infty}_{\pm},q_{\pm}^{\infty})$ are 
determined by the Milnor index
 $\tau \in {\rm Mil}_{r,s}(F)$.


\begin{prop}\label{real Hasse-Witt} We have ${\rm dim}(q_f^{\infty}) = {\rm deg}(f)n_f$, ${\rm det}(q_f^{\infty}) = [f(1)f(-1)]^{n_f}$, and
the Hasse-Witt invariant of $q_f^{\infty}$ satisfies
$$w_2(q_f^{\infty}) + w_2(Q_f) = {\rm cor}_{(E_f)_0^{\infty}/{\bf R}}({\rm det}(h^{\infty}_f),d_f)$$ in ${\rm Br}_2({\bf R})$. 

\end{prop}

\noindent
{\bf Proof.} The assertion concerning the dimension is clear, the one on the determinant follows from 
Lemma \ref{determinant}, and the property of the Hasse-Witt invariants from \cite {B 20}, Proposition 12.8.

\medskip
 For $f \in I_1$, set
$E_f^{\infty} = E_f \otimes_{\bf Q} {\bf R}$ and $(E_f)_0^{\infty} = (E_f)_0 \otimes_{\bf Q} {\bf R}$. There exists a unique non-degenerate hermitian form $(M_f^{\infty},h_f)$ over
$(E_f^{\infty},\sigma_f)$ such that $$q_f^{\infty}(x,y) = {\rm Tr}_{E_f^{\infty}/{\bf R}}(h^{\infty}_f(x,y)),$$
see for instance \cite {M}, Lemma 1.1 or \cite {B 15}, Proposition 3.6.
Set 
$$\lambda^{\infty}_f = {\rm det}(h^{\infty}_f) \in (E^{\infty}_f)_0^{\times}/{\rm N}_{E^{\infty}_f/(E^{\infty}_f)_0}.$$

\section {Obstruction group}\label{obstruction section}

We keep the notation of the previous sections. The aim of this section is to define a finite elementary abelian $2$-group
that will play an important role in the Hasse principle (see \S \ref {HP}). Recall that $J$ is the set of irreducible factors
of the polynomial $F$, that $I_0 \subset J$ is the set of factors of type 0, $I_1 \subset J$ is the set of factors of type $1$, and $I = I_0 \cup I_1$.

\begin{notation} If $f \in {\bf Z}[X]$ is an irreducible, symmetric polynomial of even degree, set $E_f = {\bf Q}[X]/(f)$, let $\sigma_f : E_f \to E_f$ be the involution
induced by $X \mapsto X^{-1}$, and let $(E_f)_0$ be the fixed field of $\sigma$ in $E_f$. Let $\alpha \in E_f$ be the image of $X$.

\end{notation}

\begin{defn}\label{ramified} Let $f \in {\bf Z}[X]$ be an irreducible, symmetric polynomial of even degree, and let $p$ be a prime number. We say that
$f$ is {\it ramified at} $p$ if there exists a place $w$ of $(E_f)_0$ above $p$ that is ramified in 
$E_f$; otherwise, we say that $f$ is {\it unramified} at $p$. We denote by $\Pi^r_f$ the set of prime numbers $p$ such that $f$ is ramified at $p$.

\medskip Let $p \in \Pi^r_f$ be an odd prime number.  If $w$ is a place of $(E_f)_0$ above $p$ that is ramified in $E_f$, we denote by $\kappa_w$  the
residue field of $w$, and by $\overline \alpha$ be the image of $\alpha$ in $\kappa_w$; we denote by $S_+$ the set of places $w$ above $p$ such that $\overline \alpha = 1$,
and by $S_-$ the set of places $w$ above $p$ such that $\overline \alpha = - 1$. We denote by 
 $\Pi^{r,+}_f$ the set of prime numbers $p$ such that there exists a place $w$ above $p$ with $w \in S_+$, and by $\Pi^{r,-}_f$ the set of 
 prime numbers $p$ such that there exists a place $w$ above $p$ with $w \in S_-$.

\end{defn}

\begin{notation} If $f, g \in {\bf Z}[X]$ are monic, irreducible, symmetric polynomials of even degree, we denote by  $\Pi_{f,g}$  the set of prime numbers $p$ such that one of the
following conditions holds :

\medskip
(a) The polynomial $f$ has a symmetric, irreducible factor $f' \in {\bf Z}_p[X]$, the polynomial $g$ has a symmetric, irreducible factor $g' \in {\bf Z}_p[X]$, such that
 $f'  \ {\rm (mod \ {\it p})}$ and $g'  \ {\rm (mod \ {\it p})}$
 have a common irreducible, 
symmetric factor  in ${\bf F}_p[X]$. 

\medskip
(b) $p \in \Pi^r_f \cap \Pi^r_g$, and the polynomials 
 $f  \ {\rm (mod \ {\it p})}$ and $g  \ {\rm (mod \ {\it p})}$
 are both divisible by $X-1$ 
 in ${\bf F}_p[X]$. 

\medskip
(c) $p \in \Pi^r_f \cap \Pi^r_g$, and the polynomials 
 $f  \ {\rm (mod \ {\it p})}$ and $g  \ {\rm (mod \ {\it p})}$
 are both divisible by $X+1$ 
 in ${\bf F}_p[X]$. 

\end{notation}

\medskip
\centerline {$F_1 = \underset{f \in I_1} \prod f^{n_f}$ and
$F_0(X) = (X-1)^{n^+} (X+1)^{n^-}$} 
\noindent for some integers $n^+,n^- \geqslant 0$.

\medskip

For all prime numbers $p$, let $D^p_+, D^p_- \in {\bf Q}_p^{\times}/{\bf Q}_p^{\times 2}$.

\begin{notation} If  $f \in I_1$, let $\Pi_{f,X-1}$ be the set of prime numbers $p$
such that $p \in \Pi^r_f$,  that $f  \ {\rm (mod \ {\it p})}$
is divisible by $X-1$ 
 in ${\bf F}_p[X]$, and that if $n^+ = 2$, then $D^p_+ \not = -1$.
 
 \medskip
 Let $\Pi_{f,X+1}$ be the set of prime numbers $p$
such that $p \in \Pi^r_f$, and that $f  \ {\rm (mod \ {\it p})}$
is divisible by $X+1$ 
 in ${\bf F}_p[X]$, and that if $n^- = 2$, then $D^p_- \not = -1$.
 
 \medskip
 Let $\Pi_{X-1,X+1} = \{2\}$ if the following conditions hold : $n^+ \not = 0$, $n^- \not = 0$, and  if $n^+ = 2$, then  $D^2_+ \not = -1$; if $n^- = 2$, then $D^2_- \not = -1$.
 Otherwise, set $\Pi_{X-1,X+1} = \varnothing$.

\end{notation}

\medskip
We denote by $C(I)$ the set of maps $I \to {\bf Z}/2{\bf Z}$.

\begin{notation}\label{f,g}
 
 If $f, g \in I$, let $c_{f,g} \in C(I)$ be such that  

\medskip
\centerline {$c_{f,g}(f) = c_{f,g}(g) = 1$, and $c_{f,g}(h) = 0$
if $h \not = f,g$. }

\medskip
Let $(f,g): C(I) \to C(I)$
be the map map sending $c$ to $c + c_{f,g}$. 

\end{notation}

\begin{notation}

Let $C_0 (I)$ be the set of $c \in C(I)$ such that

\medskip 
\centerline { $c(f) = c(g)$ \  if  \ $\Pi_{f,g} \not =  \varnothing$,} 

\medskip 
\noindent
and we denote by $\sha_F(D_+,D_-)$  the quotient of the group
 $C_0(I)$ by the subgroup of constant maps. 
 
 \medskip 
 In general,  the group  depends on $D_+ = (D_+^p)$ and $D_- = (D_-^p)$. If $n^+ \not = 2$ and $n^- \not = 2$, then $\sha_F(D_+,D_-)$ only depends on $F$, and
 we denote it by $\sha_F$.

 \end{notation}

\section{Local data}\label{local data section}

We keep the notation of \S \ref {decomposition section}.
Assume that conditions (C 1) and (C 2) of \S \ref{local-global problem section} hold, and recall that this implies the
existence of a ``local solution" everywhere. This leads, for all prime numbers $p$, to an orthogonal decomposition of the associated ${\bf Q}_p[\Gamma]$-bilinear form (see \S \ref{decomposition section}). We obtain in this way a collection of ${\bf Q}_p[\Gamma]$-bilinear forms,
one for each irreducible, symmetric factor of the characteristic polynomial. The dimensions and determinants of
the bilinear forms are always the same, but their Hasse-Witt invariants vary. 

\medskip The aim of this section is to give a combinatorial encoding of the possible Hasse-Witt invariants, called ``local data".

\medskip We identify ${\rm Br}_2({\bf R})$ and ${\rm Br}_2({\bf Q}_p)$, where $p$ is a prime number, with $\{0,1\} = {\bf Z}/2{\bf Z}$. Let $\mathcal V$ be the set of all places of ${\bf Q}$, and let $\mathcal V'$ be the set of finite places. 

\medskip
If $p$ is a prime number, let $q^p_f$ for $f \in I_1$, $q^p_+$ and $q^p_-$ be as in \S \ref{decomposition section}; recall
that if $n_+ \not = 0$ and $n_- \not = 0$, we choose $q^p_{\pm}$ such that  ${\rm det}(q_{\pm}^p) = (-1)^{s_{\pm}}|F_1(\pm 1)|$
(see Proposition \ref {dim det w pm} (b)). 

\medskip
Let $a^p \in C(I)$ be the map defined as follows :

\medskip
 $$a^p(f) = w_2(q^p_f) + w_2(Q_f)$$ if $f \in I_1$, set $$a^p(X \pm 1) = w_2(q^p_{\pm}) + w_2(Q_{\pm}).$$ 


\bigskip
Let $\mathcal C^p$ be the set of  maps $a^p \in C(I)$ obtained in this way.

 \begin{prop}\label{almost zero} For almost all prime numbers $p$, the zero map belongs to the set $\mathcal C^p$. 
 
 \end{prop}
 
 \noindent
 {\bf Proof.} Let $S$ be the set of prime numbers such that $p$ is ramified in the extension $E_f/{\bf Q}$ for some
 $f \in I_1$, or $w_2(q) \not = w_2(Q)$ in ${\rm Br}_2({\bf Q}_p)$; this is a finite set. We claim that if $p \not \in S$, then 
 the zero map belongs to $\mathcal C^p$. Indeed, set $q_f^p = Q_f^p$ for all $f \in I_1$, and 
 $q_{\pm}^p = Q_{\pm}^p$. We have ${\rm det}(q) = {\rm det}(Q)$ in $\bf Q^{\times}/{\bf Q}^{\times 2}$, and if $p \not \in S$ we have $w_2(q)  = w_2(Q)$ in ${\rm Br}_2({\bf Q}_p)$, therefore, for $p \not \in S$, we have
 
 $$(M,q) \otimes_{\bf Q} {\bf Q}_p = \underset {f \in I_1} \oplus (M^p_f,q_f^p) \oplus (M^p_+,q_+^p) \oplus (M^p_-,q_-^p) \oplus (M^p_2,q_2^p).$$
 
 If $p$ is unramified in all the extensions  $E_f/{\bf Q}$ for $f \in I_1$, by \cite{B 20}, Lemma 11.2 we
 have $\partial_p[M^p_f,q_f^p] = 0$ in $W_{\Gamma}({\bf F}_p)$; moreover,  $v_p(D_{\pm}) = 0$, hence 
 $\partial [M^p_{\pm},q_{\pm}^p] = 0$ in $W_{\Gamma}({\bf F}_p)$. 
 
 \medskip The above arguments show that if $p \not \in S$, then the choice of $q_f^p = Q_f^p$ for all $f \in I_1$ and 
 $q_{\pm}^p = Q_{\pm}^p$ gives rise to the element $a^p = 0$ of  $\mathcal C^p$; therefore the zero map is in $\mathcal C^p$,
 as claimed. This completes the proof of the proposition.

\bigskip
Proposition \ref{dim det w f} implies that if $f \in I_1$, then $a^p(f)$ is determined by ${\rm det}(h^p_f)$. Set
$\lambda^p_f = {\rm det}(h^p_f) \in (E^p_f)_0^{\times}/{\rm N}_{E^p_f/(E^p_f)_0}$. Set $E^p = \underset{f \in I_1} \prod E^p_f$ and
$E^p_0 = \underset{f \in I_1} \prod (E^p_f)_0$; the map $a^p$ is determined by
$\lambda^p  \in (E^p_0)^{\times}/{\rm N}_{E^p/E^p_0}({E^p})^{\times})$, and the quadratic forms $q^p_{\pm}$. 

\begin{notation}\label{[]} With the above notation, we set
$a^p = a^p[\lambda^p,q^p_{\pm}] = a[\lambda^p,q_+^p,q_-^p]$. 

\end{notation}

\begin{notation} If $f, g \in I$, let $c_{f,g} \in C(I)$ be such that  

\medskip 
\centerline {$c_{f,g}(f) = c_{f,g}(g) = 1$ and $c_{f,g}(h) = 0$
if $h \not = f,g$. }

\medskip Let $(f,g): C(I) \to C(I)$
be the map map sending $c$ to $c + c_{f,g}$. 

\end{notation}

\bigskip 
Recall that for all $f, g \in I$, the set $\Pi_{f,g}$ consists of the prime numbers $p$ such that $f  \ {\rm (mod \ p)}$ and $g  \ {\rm (mod \ p)}$ have a common
symmetric factor in ${\bf F}_p[X]$. 

\bigskip If $p$ is a prime number, let us consider the equivalence relation on $C(I)$ generated by the elementary equivalence

$$a \sim b \iff  b = (f,g) a \ {\rm with} \ p \in \Pi_{f,g}.$$

\medskip
 We denote by $\sim_p$ this equivalence relation. 
 
 \begin{prop}\label{equivalence class} The set  $\mathcal C^p$  is a $\sim_p$-equivalence class of $C(I)$. 
 
 \end{prop}

 \noindent
 {\bf Proof.} Set $A^p = w_2(q) + w_2(Q)$
in ${\rm Br}_2({\bf Q}_p) = {\bf Z}/2{\bf Z}$, and note that for all $a^p \in \mathcal C^p$, we have 
 $\underset{f \in J}\sum a^p(f) = A^p$.
 
 \medskip We start by proving that the set $\mathcal C^p$ is stable by the maps $(f,g)$ for $p \in \Pi_{f,g}$.
  Let $a^p[\lambda^p,q_{\pm}^p] \in \mathcal C^p$, let $f,g \in J$ be such that $p \in \Pi_{f,g}$, and let
 us show that $(f,g)(a^p[\lambda^p,q_{\pm}^p])  \in \mathcal C^p$. Note that if $f \in I_1$, then $p \in \Pi_{f,g}$ implies that  
 $(E^p_f)_0^{\times}/{\rm N}_{E^p_f/(E^p_f)_0}(E^p_f) \not = 0$. 
 Assume first that $f,g \in I$. There exist $\mu_f, \mu_g \in (E^p_f)_0^{\times}/{\rm N}_{E^p_f/(E^p_f)_0}(E^p_f)$
 such that ${\rm cor}_{(E_f)_0^p/{\bf Q}_p}(\mu_f,d_f) \not = {\rm cor}_{(E_f)_0^p/{\bf Q}_p}(\lambda_f,d_f)$
 and ${\rm cor}_{(E_f)_0^p/{\bf Q}_p}(\mu_g,d_g) \not = {\rm cor}_{(E_f)_0^p/{\bf Q}_p}(\lambda_g,d_g)$.
 Let $\mu^p \in E_0^p$ be obtained by replacing  $\lambda^p_f$ by $\mu^p_f$, $\lambda^p_f$ by $\mu^p_f$, and
 leaving the other components unchanged. 
 We have $a^p[\mu^p,q_{\pm}^p] = (f,g)(a^p[\lambda^p,q^p_{\pm}])$. Using the arguments of  \cite {B 20}, Propositions 16.5 and 22.1
 we see that  $a^p[\mu^p,q_{\pm}^p] \in \mathcal C^p$. Assume now that $f \in I_1$ and $g = X -1$. In this case, 
 the hypothesis  $p \in \Pi_{f,g}$ implies that there exists a place $w$ of $(E_f)_0$ above $p$ that
 ramifies in $E_f$, and such that the $w$-component $\lambda^w$ of $\lambda^p$ is such that with the notation
 of \cite {B 20}, \S 22, $\partial_p(\lambda^w)$ is in $W_{\Gamma}({\bf F}_p,N_+)$. 
 We modify the $w$-component of $\lambda^p$ to obtain $\mu_f  \in (E^p_f)_0^{\times}/{\rm N}_{E^p_f/(E^p_f)_0}(E^p_f)$
 such that 
 ${\rm cor}_{(E_f)_0^p/{\bf Q}_p}(\mu_f,d_f) \not = {\rm cor}_{(E_f)_0^p/{\bf Q}_p}(\lambda_f,d_f)$, and let
 $b^p$ be a quadratic form over $\bf Q_p$ with ${\rm dim}(b^p) = {\dim}(q^p_+)$, ${\rm det}(b^p) = {\det}(q^p_+)$,  and
 $w_2(b^p) = w_2(q^p_+) + 1$. We have $a^p[\mu^p,b^p,q^p_-] = (f,g)(a^p[\lambda^p,q^p_{\pm}])$.
 The arguments of  \cite {B 20}, Propositions 16.5 and 22.1 show that $a^p[\mu^p,b^p,q^p_-]  \in \mathcal C^p$. 
 
 \medskip 
 Conversely, let us show that if $a^p[\lambda^p,q_{\pm}^p]$ and $a^p[\mu^p,b_{\pm}^p]$ are in $\mathcal C^p$, then 
 $a^p[\lambda^p,q_{\pm}^p] \sim_p a^p[\mu^p,b_{\pm}^p]$.
 Let $J'$ be the set of $f \in J$ such that   $a^p[\lambda^p,q_{\pm}^p](f) \not = a^p[\mu^p,b_{\pm}^p](f)$. Since
 $\underset{h \in J}\sum a^p(h) = A^p$ for
 all $a^p \in \mathcal C^p$, the set $J'$ has an even number of elements. 
 
 \medskip
 Assume first that $p \not = 2$. This implies that for all $f \in J'$, we have $\partial_p(\lambda^p_f) \not = \partial_p(\mu^p_f)$
 and that if $f(X) = X \pm 1$, then $\partial_p(q_{\pm}^p) \not = \partial_p(b_{\pm}^p)$. Hence there exist $f, g \in J'$
 with $f \not = g$ such that $\partial_p(W_{\Gamma}({\bf Q}_p,M_f^p))$ and $\partial_p(W_{\Gamma}({\bf Q}_p,M_g^p))$ have a
  non-zero intersection. This implies that $p \in \Pi_{f,g}$. The element $f,g)(a^p[\lambda^p,q^p_{\pm}])$ differs
  from $a^p[\mu^p,b_{\pm}^p]$ in less elements than 
$a^p[\lambda^p,q_{\pm}^p]$. Since $J'$ is a finite set, continuing this way we see that $a^p[\lambda^p,q_{\pm}^p] \sim_p a^p[\mu^p,b_{\pm}^p]$.

\medskip
Suppose now that $p = 2$. Let $J''$ be the set of $f \in J'$ such that $\partial_2(\lambda^2_f) \not = \partial_2(\mu^2_f)$,
and note that $J''$ has an even number of elements. The same argument as in the case $p \not = 2$ shows
that applying maps $(f,g)$, we can assume that $J'' = \varnothing$. If $f \in J'$ and $f \not \in J''$, then
$\partial_2(\lambda^2_f)$ belongs to $W_{\Gamma}({\bf F}_2,1) \subset W_{\Gamma}({\bf F}_2)$.  Therefore $f,g  \in J'$ and $f,g  \not \in J''$, then 
$2 \in \Pi_{f,g}$. The number of these elements is also even, hence after a finite number of elementary
equivalences we see that $a^p[\lambda^p,q_{\pm}^p] \sim_p a^p[\mu^p,b_{\pm}^p]$. This completes
the proof of the proposition.

\medskip

 \begin{notation}\label {epsilon p} Let $a^p \in \mathcal C^p$, and let $c \in C(I)$. Set $$\epsilon_{a^p}(c) = 
 \underset{f \in I} \sum c(f) a^p(f).$$
 
 \end{notation}
 
 Recall from \S \ref{obstruction section} that $C_0(I)$ is the set of $c \in C(I)$ such that 
 
 \medskip \centerline {$c(f) = c(g)$ if $\Pi_{f,g} \not = \varnothing$.}
 
 \begin{lemma} Let $a^p$, $b^p$ be two elements of $\mathcal C^p$, and let $c \in C_0(I)$. Then
  $$\epsilon_{a^p}(c) = \epsilon_{b^p}(c).$$
 
 \end{lemma}

\noindent
{\bf Proof.} By Proposition \ref{equivalence class}, we have $a^p \sim_p b^p $; we can assume that $b^p = (f,g)a^p$ with 
$p \in \Pi_{f,g}.$ By definition, we have $b^p(h) = a^p(h)$ if $h \not = f,g$, $b^p(f) = a^p(f) + 1$ and 
$b^p(g) = a^p(g) + 1$. Since $c \in C_0(I)$ and $\Pi_{f,g} \not = \varnothing$, we have $c(f) = c(g)$, and this shows that
$\epsilon_{a^p}(c) = \epsilon_{b^p}(c)$, as claimed. 

\medskip 
Since $\epsilon_{a^p}(c)$ does not depend on the choice of  $a^p \in \mathcal C^p$, we set
$\epsilon^p (c) = \epsilon_{a^p}(c)$ for some $a^p \in \mathcal C^p$, and obtain a map $$\epsilon^p : C_0(I) \to {\bf Z}/2{\bf Z}.$$

\bigskip
By Proposition \ref{almost zero}, we have $\epsilon^p = 0$ for almost all prime numbers $p$.

\bigskip Let $\epsilon^{\rm finite} = \sum \epsilon^p$, where the sum is taken over all the prime numbers $p$; this is
a finite sum. Note that if $(X-1)(X+1)$ does not divide $F$, then $\epsilon^{\rm finite}$ only depends on $F$; it does
not depend of the choice of the Milnor signature $\tau$. We have a homomorphism 

$$\epsilon^{\rm finite} : C_0(F) \to {\bf Z}/2{\bf Z}.$$

\bigskip
Let $v_{\infty} \in \mathcal V$ be the unique infinite place. Recall that the forms $q^{\infty}_f$ and 
$q^{\infty}_{\pm}$  are uniquely determined by the choice of the Milnor index $\tau \in {\rm Mil}_{r,s}(F)$.
Let $a^{\infty} \in C(I)$ be the map defined as follows :

\medskip
 $$a^{\infty}(f) = w_2(q^{\infty}_f) + w_2(Q_f)$$ if $f \in I_1$,  $$a^{\infty}(X \pm 1) = w_2(q^{\infty}_{\pm}) + w_2(Q_{\pm}),$$ and 

\medskip \centerline {$a^{\infty}(f) = 0$ if $f \in J$ with $f \not \in I$, $f \not = X \pm 1$. }

\bigskip
We obtain a map

$$\epsilon^{\infty}_{\tau} :  C(I) \to {\bf Z}/2{\bf Z}$$ by setting $$\epsilon^{\infty}_{\tau}(c) =  \underset{f \in J} \sum c(f) a^{\infty}(f).$$

For $v \in \mathcal V$, set $\epsilon^v = \epsilon^p$ if $v = v_p$, and $\epsilon^v = \epsilon^{\infty}_{\tau}$ if $v = v_{\infty}$. Set

$$\epsilon_{\tau}(c) = \underset{v \in \mathcal V} \sum \epsilon^v(c).$$

Since $\epsilon^v = 0$ for almost all $v \in \mathcal V$ (cf. Proposition \ref{almost zero}), this is a finite sum. We have
$\epsilon_{\tau} = \epsilon^{\rm finite} + \epsilon^{\infty}_{\tau}$. We obtain
a homomorphism

$$\epsilon _{\tau}: C_0(I) \to {\bf Z}/2{\bf Z}.$$

\bigskip
Recall from \S \ref{obstruction section} that  $\sha_F(D_+,D_-)$   is the
quotient of $C_0(I)$ by the constant maps. 

\begin{prop}\label{epsilon} The homomorphism $\epsilon_{\tau} : C_0(I) \to {\bf Z}/2{\bf Z}$ induces a homomorphism
 $$\epsilon_{\tau} : \sha_F(D_+,D_-)  \to {\bf Z}/2{\bf Z}.$$

\end{prop}

\noindent
{\bf Proof.} It suffices to show that  if $c(f) = 1$ for all $f \in J$, then $\epsilon(c) = 0$. 
For all $v \in \mathcal V$, set $A^v = w_2(q) + w_2(Q)$
in ${\rm Br}_2({\bf Q}_v) = {\bf Z}/2{\bf Z}$, where 
${\bf Q}_v$ is either $\bf R$ or ${\bf Q}_p$, for a prime number $p$. Note that 
$A^v = 0$ for almost all $v \in \mathcal V$, and that
$\underset {v \in \mathcal V} \sum A^v = 0$. Moreover, for all $a^v \in \mathcal C^v$, we have by 
definition $\underset{f \in J}\sum a^v(f) = A^v$.

\medskip Let $c \in C(I)$ be such that $c(f) = 1$ for all $f \in J$. We have 
$$\epsilon_{\tau}(c) = \underset{v \in \mathcal V} \sum \  \underset {f \in J} \sum c(f) \ a^v(f) = \underset{v \in \mathcal V} \sum \  \underset {f \in J} \sum  \ a^v(f) = \underset {v \in \mathcal V} \sum A^v = 0.$$

\section{Hasse Principle}\label{HP}

We keep the notation of the previous sections;  in particular, $F \in {\bf Z}[X]$ is a monic, symmetric polynomial of even degree such that $F(0) = 1$ and we set $2n = {\rm deg}(F)$. 
Let $r, s \ge 0$ be integers such that $r + s = {\rm deg}(F)$ and that  $r  \equiv s \ {\rm (mod \ 8)}$, and let $\tau \in {\rm Mil}_{r,s}(F)$ be a Milnor index. 
 We assume that conditons
(C 1) and (C 2) hold.  

\medskip
Recall from \S \ref{local data section} that we have a homomorphism 
$$\epsilon_{\tau} : \sha_F(D_+,D_-)  \to {\bf Z}/2{\bf Z}.$$
 
 \begin{theo}\label{final} There exists an even, unimodular lattice having a semi-simple isometry with characteristic
 polynomial $F$ and Milnor index $\tau$ if and only if $\epsilon_{\tau} = 0$. 
 
 \end{theo}
 
 \noindent
 {\bf Proof.} Assume that there  exists an even, unimodular lattice $(L,q)$ having a semi-simple isometry with characteristic
 polynomial $F$ and Milnor index $\tau$, and let $(M,q)$ be the associated ${\bf Q}[\Gamma]$-quadratic form. Let $M^0 \oplus M^1 \oplus M^2$ the corresponding orthogonal decomposition
 of \S \ref{Q section}. We have the further orthogonal decompositions
  $(M^1,q^1) = \underset {f \in I_1} \oplus (M_f,q_f)$, and $(M^0,q^0) = (M^+,q^+) \oplus (M^-,q^-)$ (see \S 
  \ref{isometries section}). For
all prime numbers $p$, this gives rise to a local decomposition as in \S \ref{decomposition section}, and to an
element  $a^p \in \mathcal C^p$ given by $a^p(f) = w_2(q^p_f) + w_2(Q_f)$ if $f \in I_1$, by $a^p(X \pm 1) = w_2(q^p_{\pm}) + w_2(Q_{\pm}),$ and  $a^p(f) = 0$ if $f \in J$ with $f \not \in I$, $f \not = X \pm 1$ (see \S \ref{local data section}). 
Similarly, we have the element $a^{\infty} \in C(I)$ given by
 $a^{\infty}(f) = w_2(q^{\infty}_f) + w_2(Q_f)$ if $f \in I_1$,  $a^{\infty}(X \pm 1) = w_2(q^{\infty}_{\pm}) + w_2(Q_{\pm}),$ and 
$a^{\infty}(f) = 0$ if $f \in J$ with $f \not \in I$, $f \not = X \pm 1$. 
Since $q^p_f = q_f \otimes_{\bf Q}{\bf Q}_p$ and $q^{\infty}_f = q_f \otimes_{\bf Q}{\bf R}$
for all $f \in J$, we have

\medskip

\centerline { $ \underset{v \in \mathcal V} \sum \  a^v(f) = 0$ for all $f \in J$.}

\medskip This implies that $\epsilon_{\tau} = 0$.

\bigskip Conversely, assume that $\epsilon_{\tau} = 0$.
By \cite{B 20}, Theorem 13.5 this implies that for all $v \in \mathcal V$ there exists $b^v \in \mathcal C^v$ such
that 
for all $f \in J$, we have $\underset {v \in \mathcal V} \sum b^v(f) = 0$. 

\medskip
If $v \in \mathcal V'$ with $v = v_p$ where  $p$ is a prime number, let us write $b^v = a[\lambda^p,q_+^p,q_-^p]$, for some
$\lambda^p  \in (E^p_0)^{\times}/{\rm N}_{E^p/E^p_0}({E^p})^{\times}$, and some  quadratic forms 
$q_+^p,q_-^p$ over ${\bf Q}_p$, as in notation \ref{[]}. 

\medskip
Note that since $v_{\infty}$ does not belong to any of the sets $\Pi_{f,g}$, 
we have $b^{\infty}(f) = a^{\infty}(f) = w_2(q^{\infty}_f) + w_2(Q_f)$ if $f \in I_1$,  $b^{\infty}(f) = a^{\infty}(X \pm 1) = w_2(q^{\infty}_{\pm}) + w_2(Q_{\pm}),$ and 
$b^{\infty}(f) = a^{\infty}(f) = 0$ if $f \in J$ with $f \not \in I$, $f \not = X \pm 1$.
Recall that the forms $q^{\infty}_f$ and 
$q^{\infty}_{\pm}$  are uniquely determined by the choice of the Milnor index $\tau \in {\rm Mil}_{r,s}(F)$.

 \bigskip
 If $f \in I_1$, we have $$b^{v_p}(f) = a[\lambda^p,q_+^p,q_-^p](f) = {\rm cor}_{(E_f)_0^p/{\bf Q}_p}(\lambda^p_f,d_f),$$ and 
 $$b^{v_{\infty}}(f) = a[\lambda^{\infty},q_+^{\infty},q_-^{\infty}](f) =  {\rm cor}_{(E_f)_0^{\infty}/{\bf R}}({\rm det}(h^{\infty}_f),d_f)$$
 
 \bigskip
 \noindent
  (cf. Propositions \ref{dim det w f} and  \ref{real Hasse-Witt}). 
  
  \medskip Since $\underset {v \in \mathcal V} \sum b^v(f) = 0$, we have 
  $$\underset {v \in \mathcal V} \sum {\rm cor}_{(E_f)_0^v/{\bf Q}_v}(\lambda^v_f,d_f) = 0,$$
  where ${\bf Q}_v = {\bf Q}_p$ if $v = v_p$ and ${\bf Q}_v = {\bf R}$ if $v = v_{\infty}$. This implies that
   $$\underset {w \in {\mathcal W}} \sum (\lambda^w_f,d_f) = 0,$$ where $\mathcal W$ is the set of primes of $E_0$. 
   Therefore there exists $\lambda_f \in E_0^{\times}/N_{E/E_0}(E^{\times})$ mapping to $\lambda_f^w$ for all
   $w \in {\mathcal W}$ (see for instance \cite{B 20}, Theorem 10.1). In particular, we have $(\lambda_f,d_f) = (\lambda^w_f,d_f)$
   in ${\rm Br}_2(E_0^w)$ for all $w \in {\mathcal W}$. 
   
   \medskip
   Note that $\tau(f)$ is an even integer. Let $h_f : M_f \times M_f \to E_f$ be a hermitian form such that
   ${\rm det}(h_f) = \lambda_f$, and that the index of $h_f$ is equal to ${\tau (f)} \over 2$; such a hermitian
   form exists (see for instance  \cite{Sch}, 10.6.9). Let us define   $$q_f : M_f \times M_f \to {\bf Q}$$ by $$q_f(x,y) = {\rm Tr}_{E_f/{\bf Q}}(h_f(x,y)).$$

   \bigskip
   Let $f = X \pm 1$. We have  $\underset {v \in \mathcal V} \sum b^v(f) = 0$,  hence by the Brauer-Hasse-Noether theorem
   there exists $a(\pm) \in {\rm Br}_2({\bf Q})$ mapping to $b^v(f)$ in ${\rm Br}_2({\bf Q}_v)$ for all $v \in \mathcal V$.
   Let $q_{\pm}$ be a quadratic form over $\bf Q$ of dimension $n_{\pm}$, determinant $D_{\pm}$, Hasse-Witt
   invariant $w_2(q_{\pm}) = a(\pm) + w_2(Q_{\pm})$ and index $\tau(X\pm 1) = r_{\pm}-s_{\pm}$. 
   Such a quadratic form exists; see for instance \cite{S}, Proposition 7.
   
   \medskip
   Let $q' : M \times M \to {\bf Q}$ be the quadratic form given by $$(M,q')= \underset{f \in I_1} \oplus (M_f,q_f) \oplus 
   \underset {f \in I_0} \oplus (M_f,q_f) \oplus (M^2,q^2),$$ where $(M^2,q^2)$ is hyperbolic. By construction,
   $(M,q')$ has the same dimension, determinant, Hasse-Witt invariant and signature as $(M,q)$, hence
   the quadratic forms $(M,q')$ and $(M,q)$ are isomorphic. 
   
   \medskip Let $t : M \to M$ be defined by $t(m) = \gamma m$, where $\gamma$ is a generator of $\Gamma$. 
  By construction, $t$ is an isometry of $(M,q')$ and it is semi-simple with characteristic 
   polynomial $F$. 
   By hypothesis, conditions (C 1) and (C 2) hold, hence $(M,q')\otimes_{\bf Q} {\bf Q}_p$ contains
   an even, unimodular ${\bf Z}_p$ lattice $L_p$ stable by the isometry $t$. Let
   $$L = \{ x \in M \ | \ x \in L_p \ {\rm for \ all} \ {\rm prime \ numbers} \ p \}.$$  
 $(L,q')$ is an even, unimodular lattice having a semi-simple isometry with characteristic polynomial $F$. This
 completes the proof of the theorem.
 
 \begin{coro}\label{final coro} Assume that conditions {\rm (C 1)} and {\rm (C 2)} hold, and that $\sha_F = 0$. Then there exists
 an even, unimodular lattice having a semi-simple isometry with characteristic polynomial $F$ and Milnor index $\tau$. 
 
 \end{coro}

\section{Even, unimodular lattices preserved by a semi-simple element of ${\rm SO}_{r,s}({\bf R})$ }\label{reformulation}

In this section, we reformulate the Hasse principle result of \S \ref{HP}, and prove a result stated in the introduction.
We keep the notation of \S \ref{HP}. In particular $F \in {\bf Z}[X]$ is a monic, symmetric polynomial of even degree such that $F(0) = 1$, and
$r, s \ge 0$ are integers such that $r + s = {\rm deg}(F)$ and that  $r  \equiv s \ {\rm (mod \ 8)}$.

\medskip Let us now assume that condition (C 2) holds, and let $t \in {\rm SO}_{r,s}({\bf R})$ be a semi-simple
 isometry with characteristic polynomial $F$. Let $\tau = \tau(t) \in {\rm Mil}_{r,s}(F)$ be the Milnor index
 associated to $t$ in Proposition \ref{bijection}. 
 
 \medskip  Assume  that condition (C 1) also holds, and let $\epsilon_{\tau} : \sha_F(D_+,D_-)  \to {\bf Z}/2{\bf Z}$ be the
 homomorphism defined in \S \ref{HP}; set $\epsilon_t = \epsilon_{\tau}$. 
 The following is a reformulation of Theorem \ref{final} :
  
 \begin{theo}\label{preserves} The isometry  $t \in {\rm SO}_{r,s}({\bf R})$ preserves an even, unimodular lattice if and
 only if $\epsilon_t = 0$.
 
 \end{theo}
 
  \begin{coro}\label{preserves coro} If $\epsilon_{\tau} : \sha_F(D_+,D_-)  \to {\bf Z}/2{\bf Z} = $,  the isometry  $t \in {\rm SO}_{r,s}({\bf R})$ preserves an even, unimodular lattice.
 
 \end{coro}

\section{Automorphisms of $K3$ surfaces}\label{K3}

Which Salem numbers occur as dynamical degrees of automorphisms of complex analytic $K3$ surfaces ? This question was raised
by Curt McMullen in \cite{Mc1}, and was studied in many other papers (see for instance \cite {GM}, \cite {O}, \cite {Mc2}, \cite{R},
\cite {BG}, \cite {Mc3}, \cite{R2}, \cite{Z},  \cite {Br}). 

\medskip
We refer to \cite{H} and  \cite{Ca} for background on complex $K3$ surfaces (henceforth
$K3$ surfaces, for short) and their automorphisms. 

\medskip
Let $\mathcal X$ be a $K3$ surface, and let $T : \mathcal X \to \mathcal X$ be an automorphism; it induces an isomorphism
$T^* : H^2(\mathcal X,{\bf Z}) \to H^2(\mathcal X,{\bf Z})$. The {\it dynamical degree} of $T$ is by definition 
the spectral radius of $T^*$; it is either 1 or a Salem number. The characteristic polynomial
of $T^*$ is a product of at most one Salem polynomial and of a finite number of cyclotomic polynomials (see 
\cite {Mc1}, Theorem 3.2). 

\medskip Let $H^2(\mathcal X,{\bf C}) = H^{2,0}(\mathcal X) \oplus H^{1,1}(\mathcal X) \oplus H^{0,2}(\mathcal X)$ be the Hodge decomposition of $H^2(\mathcal X,{\bf C})$. Since the subspace $H^{2,0}(\mathcal X)$ is one dimensional, $T^*$ acts
on it by multiplication by a scalar, denoted by $\delta(T)$, and called the {\it determinant} of $T$; we have $|\delta(T)| = 1$. Moreover, $\delta(T)$ is a root of unity if  $\mathcal X$ is projective (cf. \cite {Mc1}, Theorem 3.5). 




\medskip The intersection form of $H^2(\mathcal X,{\bf Z)}$ is an even, unimodular lattice of signature (3,19), hence it is
isomorphic to $\Lambda_{3,19}$, and an automorphism of $\mathcal X$ induces an isometry of that form. Therefore a
necessary condition for a Salem number $\alpha$  to occur as the dynamical degree of such an automorphism is that $\Lambda_{3,19}$ has an isometry with
characteristic polynomial $S C$, where $S$ is the minimal polynomial of $\alpha$, and $C$ is a (possibly empty) product
of cyclotomic polynomials.

\begin{defn}\label{complemented} A {\it complemented Salem polynomial} is by definition a degree $22$ polynomial that
is the product of a Salem polynomial and of a (possibly empty) product of cyclotomic polynomials. 

\end{defn}

Recall  from \S \ref{local-global problem section}  that a monic, symmetric polynomial $F \in {\bf Z}[X]$ satisfies condition (C~1) if and only if

\medskip 
\centerline  {\it $|F(1)|$, $|F(-1)|$ and $(-1)^n F(1) F(-1)$ are squares,}

\medskip
\noindent
where $2n = {\rm deg}(F)$, and that this condition is {\it necessary} for $F$ to be the characteristic polynomial of an isometry of an even, unimodular lattice.

\medskip
If $F$ is a complemented Salem polynomial, then $m(F) = 1$, since $F$ has exactly two roots that are not on the unit circle. This implies that condition (C 2) holds for  $(r,s) = (3,19)$.

 \begin{defn}\label{realizable polynomial} Let $F$ be a complemented Salem polynomial, and let $\delta$ be a root of $F$
with $|\delta| = 1$. We say that $(F,\delta)$ is {\it realizable} (resp. projectively realizable) if there exists
a  $K3$ surface (resp. a projective $K3$ surface)  $\mathcal X$
 and an automorphism $T : \mathcal X \to \mathcal X$ such
 that
 
 \medskip \noindent
 $\bullet$
 \ \ $F$ is the characteristic polynomial of $T^*|H^2(\mathcal X).$ 
 
  \medskip \noindent
$\bullet$ \ \ $T^*$ acts on $H^{2,0}(\mathcal X)$ by multiplication by $\delta$.

\end{defn}

Let $S$ be a Salem polynomial of degree $d$ with  $4 \leqslant d \leqslant 20$, and set  $$F(X) = S(X)(X-1)^{22-d}.$$
Let us consider Salem polynomials $S$ with
$|S(1)| = 1$. In this case, $\Pi_{S(X),X-1} = \varnothing$, hence  the obstruction group $\sha_F$ is not trivial, and not all Milnor indices are realized. We start 
by introducing some notation.

\begin{notation}
If $\delta$ is a root of $S$ with $|\delta| = 1$, let $\tau_{\delta} \in {\rm Mil}_{3,19}(F)$
be such that 

\medskip

\centerline {$\tau_{\delta}(\mathcal P) = 2$ \ \  if \ \ $\mathcal P(X) = (X-\delta)(X-\delta^{-1})$,}

\medskip
\noindent
that

\medskip

\centerline {$\tau_{\delta}(\mathcal Q) = -2$ \ \ for all \ \ $\mathcal Q \in {\rm Irr}_{\bf R}(S)$ with $\mathcal Q \not = \mathcal P$,}

\medskip
\noindent
and that 
$$\tau_{\delta}(X-1) = d - 22.$$

\bigskip

Let $\tau_{1} \in {\rm Mil}_{3,19}(F)$ be such that 
$\tau_{1}(\mathcal Q) = -2$ for all $\mathcal Q \in {\rm Irr}_{\bf R}(S)$, and that 
$\tau_{1}(X-1) = d - 20$.

\end{notation}

\begin{theo}\label{theorem 4 - second part}
  Let $S$ be a Salem polynomial of degree $d$ with $4 \leqslant d\leqslant 20$, set
 $$F(X) = S(X)(X-1)^{22-d}.$$
 Assume that condition {\rm (C 1)} holds for $F$ and that $|S(1)| = 1$.
 Let
 $\tau \in {\rm Mil}_{3,19}(F)$.
 
 \medskip
 Then the lattice $\Lambda_{3,19}$ has an isometry with characteristic polynomial $F$ and Milnor index $\tau$ if
 and only if one of the following holds

\medskip
{\rm (i)} $d \equiv -2 \ {\rm (mod \ 8)}$ and $\tau = \tau_{\delta}$ where $\delta$ is a root of $S$ with $|\delta| = 1$.

\medskip
{\rm (ii)} $d \equiv 2 \ {\rm (mod \ 8)}$ and $\tau = \tau_{1}$.

\end{theo}

\medskip
\noindent
{\bf Proof.} The polynomials $S$ and $X-1$ are relatively prime over ${\bf Z}$. This implies that if the lattice 
$\Lambda_{3,19}$ has a semi-simple isometry with characteristic polynomial $F$, then $\Lambda_{3,19} \simeq L_1 \oplus L_2$
where $L_1$ and $L_2$ are even, unimodular lattices, such that $L_1$ has an isometry with characteristic polynomial $S$, 
and $L_2$
has a semi-simple isometry with characteristic polynomial $(X-1)^{22-d}$. 


\medskip Note that every $\tau \in {\rm Mil}_{3,19}(F)$ is either equal to $\tau_1$, or to $\tau_{\delta}$, where $\delta$ is
a root of $S$ with $|\delta| = 1$.
Assume first  that $\Lambda_{3,19}$ has a semi-simple isometry with characteristic polynomial $F$ and
Milnor index $\tau_{1}$, and let $\Lambda_{3,19} \simeq L_1 \oplus L_2$ be as above. The signature of $L_1$ is  $(1,d-1)$, and since $L_1$ is unimodular and
even, this implies that  $d \equiv 2 \ {\rm (mod \ 8)}$.

\medskip Suppose now that $\Lambda_{3,19}$ has a semi-simple isometry with characteristic polynomial $F$ and
Milnor index $\tau_{\delta}$, where $\delta$ is
a root of $S$ with $|\delta| = 1$. Let
$\Lambda_{3,19} \simeq L_1 \oplus L_2$ be as above. The signature of $L_1$ is then $(3,d-3)$, and since $L_1$ is unimodular and
even, we have $d \equiv -2 \ {\rm (mod \ 8)}$.

\medskip This implies that if $\Lambda_{3,19}$ has a semi-simple isometry with characteristic polynomial $F$, then we
are in one of the cases (i) or (ii). 

\medskip
Let us show the converse. Suppose first that we are in case (i). We have $d \equiv -2 \ {\rm (mod \ 8)}$; this means that $d = 6$ or $d = 14$. 
Let $(r,s) = (3,3)$ if $d = 6$ and $(r,s) = (3,11)$ if $d = 14$; note that condition (C 2) holds for $S$ and $(r,s)$, and that $r  \equiv s \ {\rm (mod \ 8)}$. By hypothesis, condition (C 1)  holds for $F$; since $F(-1) = S(-1)$, this implies that $|S(-1)|$ is a square. 
Moreover, $|S(1)| = 1$ by hypothesis. We claim that condition (C~1) also holds for $S$. 
Since $S$ is a Salem polynomial, we have $S(1) < 0$ and $S(-1) > 0$; we have $d  \equiv 2 \ {\rm (mod \ 4)}$, therefore
$(-1)^{d/2}S(1)S(-1)$ is a square. This implies that condition (C 1) holds for $S$. Moreover, $S$ is irreducible, hence
$\sha_S = 0$.

\medskip
Let $\tau' \in {\rm Mil}_{r,s}(S)$ be the restriction of $\tau_{\delta}$ to ${\rm Mil}_{r,s}(S)$. We have seen that
conditions (C 1) and (C 2) hold for $S$, and that $\sha_S = 0$. By Corollary \ref {final coro} the even, unimodular lattice $\Lambda_{r,s}$ has an
isometry with characteristic polynomial $S$ and Milnor index $\tau'$. The identity is a semi-simple isometry of the lattice $-E_8$ with
characteristic polynomial $(X-1)^{8}$. Since  $\Lambda_{3,19} = \Lambda_{3,3} \oplus (-E_8) \oplus (- E_8)
=  \Lambda_{3,11} \oplus (-E_8)$,  the lattice $\Lambda_{3,19}$ has a semi-simple isometry with characteristic polynomial $F$ and Milnor index $\tau_{\delta}$, as claimed.

\medskip
Suppose now that we are in case (ii). We have $d \equiv 2 \ {\rm (mod \ 8)}$; that is,  $d = 10$ or $d = 18$. 
Let $(r,s) = (1,d-1)$;
note that $r  \equiv s \ {\rm (mod \ 8)}$, and that condition (C 2) holds for $S$ and $(r,s)$. We show as in case (i) that condition (C 1) holds for $S$ and that $\sha_S = 0$.

\medskip
Let $\tau'' \in {\rm Mil}_{r,s}(S)$ be the restriction of $\tau_{1}$ to ${\rm Mil}_{r,s}(S)$. By Corollary \ref {final coro} the even, unimodular lattice $\Lambda_{r,s}$ has an
isometry with characteristic polynomial $S$ and Milnor index $\tau''$. The identity is a semi-simple isometry of the lattice $\Lambda_{2,20-d}$ with
characteristic polynomial $(X-1)^{22-d}$. Since  $\Lambda_{3,19} = \Lambda_{1,d-1} \oplus \Lambda_{2,20-d}$,  the lattice $\Lambda_{3,19}$ has a semi-simple isometry with characteristic polynomial $F$ and Milnor index $\tau_{1}$.



\begin{prop}\label{summary}
 Let $S$ be a Salem polynomial of degree $d$ with $4 \leqslant d\leqslant 22$ and $d \equiv 6 \ {\rm (mod \ 8)}$.
 and let $\delta$ be a root of $S$ with $|\delta| = 1$. Suppose that $|S(1)|$ and $S(-1)$ are both squares, and set $F(X) = S(X)(X-1)^{22-d}$. Then
 $(F,\delta)$ is realizable.

\end{prop}

\noindent
{\bf Proof.} The argument of 
Theorem \ref {theorem 4 - second part} implies that 
the lattice $\Lambda_{3,19}$ has a semi-simple isometry
with characteristic polynomial $F$ and Milnor index $\tau_{\delta}$. Applying \cite {Br}, Lemma 3.3 (1) of Brandhorst gives
the desired result.

\begin{example}\label {Salem 6} If  $a \geqslant 0$ is an integer, the polynomial 
$$S_a(X) = X^6 -aX^5 - X^4 + (2a -1) X^3  - X^2 -a X + 1$$ 
is an  Salem polynomial (see \cite {GM}, page 284, Example 1), and $S(1) = -1$. Part (iii) of Theorem \ref{summary} implies
that if $\delta_a$ is a root of $S_a$ with $|\delta_a| = 1$ and $F_a(X) = S_a(X)(X-1)^{16}$, then $(F_a,\delta_a)$ is realizable. 

\medskip The  polynomials $S_a$ also appear in \S 4 of \cite {Mc1} : for every integer $a \geqslant 0$, McMullen
gives a geometric construction of an automorphism of a non-projective $K3$ surface such that the dynamical
degree and the determinant of the automorphisms are roots of $S_a$ (see \cite{Mc1}, Theorem 4.1); this construction uses
complex tori. 

\end{example}

\begin{example}\label{14} Let $S$ be a Salem polynomial of degree $14$, and assume that  $|S(1)|= 1$. Let $\delta$ be
a root of $S$, and set $F(X) = S(X)(X-1)^{8}$. If moreover $S(-1)$ is a square, then condition (C 1) holds for $F$, and Theorem \ref{summary} (iii) implies that $(F,\delta)$ is  realizable. 

\medskip
Salem polynomials of degree $14$ with  $|S(1)| = |S(-1)| = 1$ were considered in several papers. Oguiso proved that
the third smallest known Salem number $\lambda_{14}$ is the dynamical degree of an automorphism of a non-projective
$K3$ surface (see \cite {O}, Proposition 3.2). If a Salem number is a root of a Salem polynomial $S$ of degree $14$
with $|S(1)| = |S(-1)| = 1$, then this was shown by Reschke (see \cite{R}, Theorem 1.2).

\end{example}

\section{Realizable Salem numbers}\label{4}

The aim of this section is to show that if $\alpha$ is a Salem number of degree $d$ with $d = 4,6,8, 12, 14$ or $16$, then
$\alpha$ is the dynamical degree of an automorphism of a {\it non-projective} $K3$ surface;  partial results are given
for the other values of $d$ as well. 

\begin{notation} Let $S$ be a Salem polynomial of degree $d$ with $4 \leqslant d \leqslant 20$, and let $\delta$ be a root of $S$
with $|\delta| = 1$. Let $F$ be a complemented Salem polynomial with Salem factor $S$. We define the Milnor index $\tau_{\delta} \in
{\rm Mil}_{3,19}(F)$ as follows : 

\medskip
$\bullet$ $\tau_{\delta}((X - \delta)(X - \delta^{-1})) = 2$;

\medskip
$\bullet$ $\tau_{\delta}(\mathcal P)< 0$ for all $\mathcal P \in {\rm Irr}_{\bf R}(F)$ such that $\mathcal P (X) \not = 
(X - \delta)(X - \delta^{-1})$.

\end{notation}

\begin{theo}\label{n+} Let $S$ be a Salem polynomial of degree $d$, and let $\delta$ be a root of $S$
with $|\delta| = 1$. Suppose that  $d \leqslant 16$ and $d  \equiv 0, 4 \ {\rm or} \ 6  \ {\rm (mod \ 8)}$.
Then there exists a complemented Salem polynomial $F$ with Salem factor $S$ such that $\Lambda_{3,19}$ has a realizable isometry with Milnor index $\tau_{\delta}$. 

\end{theo}

\begin{notation} If $f \in {\bf Z}[X]$ is an irreducible, symmetric polynomial of even degree, set $E_f = {\bf Q}[X]/(f)$, let $\sigma_f : E_f \to E_f$ be the involution
induced by $X \mapsto X^{-1}$, and let $(E_f)_0$ be the fixed field of $\sigma$ in $E_f$.

\end{notation}

\begin{defn}\label{ramified} Let $f \in {\bf Z}[X]$ be an irreducible, symmetric polynomial of even degree, and let $p$ be a prime number. We say that
$f$ is {\it ramified at} $p$ if there exists a place $w$ of $(E_f)_0$ above $p$ that is ramified in 
$E_f$.

\end{defn}

\begin{prop}\label{not both squares} Suppose that  $d \leqslant 18$ and that one of the following holds :

\medskip {\rm (i)} $|S(1)|$ and $S(-1)$ are not both squares. 

\medskip {\rm (ii)} $S$ is ramified at the prime $2$.

\medskip
Then there exists a complemented Salem polynomial $F$ with Salem factor $S$ such that $\Lambda_{3,19}$ has a realizable isometry with Milnor index $\tau_{\delta}$. 

\end{prop}

\noindent
{\bf Proof.} (a) Suppose first  that  there exists a prime number $p$ such that $v_p(S(-1)) \equiv 1 \ {\rm (mod \ 2)}$. 
Set $F(X) = S(X)(X-1)^{20-d}(X+1)^2$. We have $D_- = (-1)^{s_-}S(-1) = S(-1)$; this implies that $D_- \not = -1$ in 
${\bf Q}_p^{\times}/{\bf Q}_p^{\times 2}$, since $v_p(D_-) = v_p(S(-1)) \equiv 1 \ {\rm (mod \ 2)}$. Therefore $p \in \Pi_{S,X+1}$, and $\sha_F(D_-) = 0$.

\medskip (b) Suppose now that no such prime number exists; this implies that $S(-1)$ is a square. Then either $|S(1)|$ is  not a square, or $S$ is ramified at $2$.

Set $F(X) = S(X)(X-1)^{22-d}$. We have $d \leqslant 18$, hence $22-d \not = 2$; this implies that $p \in \Pi_{S,X-1}$ and hence
$\sha_F = 0$. 

\medskip

Let $(L,q)$ be an even, unimodular lattice of signature $(3,19)$, and let $t : L \to L$ be a semi-simple
isometry with characteristic polynomial $F$ and with Milnor index $\tau_{\delta}$; such an isometry exists 
$\sha_F(D_-) = 0$ in case (a) and $\sha_F = 0$ in
case (b).

\medskip
Set $L_S = {\rm Ker}(S(t))$, and let $S_C$ be the orthogonal of $L_S$ in $L$. If we are in case (b),
then the restriction of $t$ to $L_C$ is the identity, hence it is a positive isometry in the terminology of McMullen; this implies that
$t : L \to L$ is realizable.

\medskip
Suppose now that we are in case (a). Let $L_1 = L_S = {\rm Ker}(S(t))$,
$L_2 = {\rm Ker}(t+1)$ and $L_3 = {\rm Ker}(t-1)$; let $t_i : L_i \to L_i$ be the restriction of $t$ to $L_i$.

\medskip A root of $(L_2,q)$ is by definition an element $x \in L_2$ such that $q(x,x) = -2$. If $(L_2,q)$  has no
roots, then $t_2$ is a positive isometry of $(L_2,- q)$ by \cite{Mc2}, Theorem 2.1, and hence \cite {Mc2}, Theorem 6.2 (see
also \cite {Mc3}, Theorem 6.1) implies that $(F_+,\delta)$ is realizable, hence  (b) holds.

\medskip Suppose that $(L_2,q)$ has at least one root.
By hypothesis, $S(-1)$ is not a square, therefore ${\rm det}(L_2,q)$ is not a square. Since $(L_2,q)$ is of rank $2$, even and negative
definite, there exist integers $D \geqslant 1$ and $f \geqslant 1$  such that 
${\rm det}(L_2,q) = f^2D$, where $-D$ is the discriminant of an imaginary quadratic field.
The lattice $(L_2,q)$  is isomorphic to a quadratic form $q'$ on an
order $O$ of the imaginary quadratic field ${\bf Q}(\sqrt {-D})$ (see for instance \cite {Co}, Theorem 7.7). Complex
conjugation induces an isometry of the quadratic form $(O,q')$ with characteristic polynomial $X^2 - 1$. If $D = 3$ and $f = 1$, then
$(O,q')$  is isomorphic to the root lattice $A_2$, and complex conjugation is a positive isometry of $(O, -q')$ (see \cite {Mc2}, \S 5, Example);
otherwise, $(O,q')$ contains only two roots,  fixed by complex conjugation, hence we obtain a positive  isometry of $(O, -q')$ in
this case as well. Let $t_2' : L_2 \to L_2$ be the isometry of $(L_2,q)$ obtained via the isomorphism $(O,q') \simeq (L_2,q)$. 
Then $t_2'$ is a positive isometry of $(L_2,-q)$. Let $G(L_2) = (L_2)^{\sharp}/L_2$, and note that $t_2$ and $t_2'$ both
induce $- {\rm id}$ on $G(L_2)$. This implies that $(L,q)$ has a semi-simple isometry $t' : L \to L$ inducing  the positive
isometry $t_2$ or $t_2'$ on
$L_2$ and $t_i$ on $L_i$ for $i = 1,3$. By  \cite {Mc2}, Theorem 6.2 (see
also \cite {Mc3}, Theorem 6.1) this implies that  the isometry is realizable.

\medskip
\noindent
{\bf Proof of Theorem \ref{n+}.} If $|S(1)|$ and $S(-1)$ are not both squares or if if $d$ is divisible by $4$, then the result follows from the proposition; if
$d \equiv 6 \ {\rm (mod \ 8)}$, then it follows from Proposition \ref{summary}.

 \begin{defn}\label{realizable}
 Let $\alpha$ be a Salem number and let $\delta$  be a conjugate of $\alpha$ such that $|\delta| = 1$. We say
that $(\alpha,\delta)$ is {\it realizable} (resp. projectively realizable) if there exists an automorphism of a  $K3$ surface
(resp. a projective $K3$ surface) 
having an automorphism of dynamical degree $\alpha$ and and determinant $\delta$.

 \end{defn}
 
 \begin{coro}\label{lots coro} Let $\alpha$ be a Salem number of degree $d$ with $4 \leqslant d \leqslant 16$, let $S$ be
 the minimal polynomial of $\alpha$, and let 
$\delta$ be a root of $S$
with $|\delta| = 1$.  If 
$d = 4, 6, 8, 12, 14$ or $16$, then
$(\alpha,\delta)$ is realizable. 
 
 \end{coro}
 
 \begin{remark}\label{20} If $d = 20$ with $|S(1)|$ is a square and $S(-1)$ is not a square, then the method of Proposition \ref{not both squares} still works,
and therefore Corollary \ref{lots coro} holds in this case as well. 

\end{remark}

\begin{example} McMullen proved that the Salem numbers $\lambda_{14}$, $\lambda_{16}$ and $\lambda_{20}$ 
are not realized as dynamical degrees of  automorphisms of {\it projective} $K3$ surfaces (cf. \cite{Mc3}, \S 9).
Corollary \ref{lots coro} and Remark \ref{20} show that they are realized by automorphisms of {\it non-projective } $K3$ surfaces. 

\end{example} 

\section{A nonrealizable Salem number}\label{18}

McMullen proved that the Salem number $\lambda_{18} = 1.1883681475...$  (the second smallest known Salem number)  is the dynamical degree of
an automorphism of a  {\it projective} $K3$ surface (cf. \cite{Mc3}, Theorem 8.1). The aim of this section is to show that this is not possible
for {\it non-projective} $K3$ surfaces.

\medskip

Let $S$ be a Salem polynomial of degree $18$, and let
$\delta$ be a root of $S$ with $|\delta| = 1$. Let $\sigma_{\delta} \in {\rm Mil}_{3,15}(S)$ 
be such that $\sigma_{\delta}(\mathcal P) = 2$ for $\mathcal P(x) = (x-\delta)(x - \delta^{-1})$ and that 
$\sigma_{\delta}(\mathcal Q) = -2$ for all $\mathcal Q \in {\rm Irr}_{\bf R}(S)$ with $\mathcal Q \not = \mathcal P$.

\medskip
If $f \in {\bf Z}[X]$ is a monic polynomial, we denote by  ${\rm Res}(S,f)$ the resultant of the polynomials $S$ and $f$. 

\begin{prop}\label{resultant} Assume that $|{\rm Res}(S,f)| = 1$ for all $f \in \{\Phi_1,\Phi_2,\Phi_3,\Phi_4, \Phi_6 \}$. Let $C$ be a product
of cyclotomic polynomials such that ${\rm deg}(C) = 4$, and set $F = SC$. Let $\tau_{\delta} \in {\rm Mil}_{3,19}(F)$ be
such that the restriction of $\tau_{\delta}$ to $ {\rm Mil}_{3,15}(S)$ is $\sigma_{\delta}$, and that 
$\tau_{\delta}(\mathcal Q) < 0$ for all $\mathcal Q \in {\rm Irr}_{\bf R}(C)$.

\medskip If $\Lambda_{3,19}$ has a semi-simple isometry with characteristic polynomial $F = SC$ and Milnor index $\tau_{\delta}$,
then $C = \Phi_{12}$.

\end{prop}

\noindent
{\bf Proof.} If $C = \Phi_5$, $\Phi_8$ or $\Phi_{10}$, then $FC$ does not satisfy Condition (C 1), hence 
$\Lambda_{3,19}$ does not have any isometry with characteristic polynomial $F$ for these choices of $C$.

\medskip
Assume that all the factors of $C$ belong to the set $\{\Phi_1,\Phi_2,\Phi_3,\Phi_4, \Phi_6 \}$. Then $S$ and $C$
are relatively prime over ${\bf Z}$. If $\Lambda_{3,19}$ has an isometry with characteristic polynomial $F$, then
$\Lambda_{3,19} = L_1 \oplus L_2$, where $L_1$ and $L_2$ are even, unimodular lattices such that
$L_1$ has an isometry with characteristic polynomial $S$ and Milnor index $\sigma_{\delta}$, and $L_2$ has
an isometry with characteristic polynomial $C$. This implies that the signature of $L_1$ is $(3,15)$ and that
the signature of $L_2$ is $(0,4)$, and this is impossible. 

\medskip Therefore the only possiblity is $C = \Phi_{12}$, as claimed. 

\begin{notation}
Let $C = \Phi_{12}$, and set $F = SC$. Let $\zeta$ be a primitive $12$th root of unity. Let $\tau_{\delta}, \tau_{\zeta}
 \in {\rm Mil}_{3,19}(F)$ be such that 
 
 \medskip
 $\tau_{\delta}(\mathcal P) = 2$ for $\mathcal P(x) = (x-\delta)(x - \delta^{-1})$ and that 
$\tau_{\delta}(\mathcal Q) = -2$ for all $\mathcal Q \in {\rm Irr}_{\bf R}(F)$ with $\mathcal Q \not = \mathcal P$;

\medskip
$\tau_{\zeta}(\mathcal P) = 2$ for $\mathcal P(x) = (x-\zeta)(x - \zeta^{-1})$ and that
$\tau_{\zeta}(\mathcal Q) = -2$ for all $\mathcal Q \in {\rm Irr}_{\bf R}(F)$ with $\mathcal Q \not = \mathcal P$.

\end{notation}

\begin{theo}\label{phi 12} Let $S$ be an  Salem polynomial of degree $18$ such that $$|S(1)S(-1)| = 1,$$ let $C = \Phi_{12}$, and set $F = SC$. 
Let $\delta$ be a root of $S$ with $|\delta| = 1$, and let $\zeta$ be a primitive $12$th root of unity. 
With the above notation, we have

\medskip

{\rm (a)} The lattice $\Lambda_{3,19}$ has an isometry with characteristic polynomial $F$ and Milnor
index $\tau_{\zeta}$.

\medskip 

{\rm (b)} The lattice $\Lambda_{3,19}$ has an isometry with characteristic polynomial $F$ and Milnor
index $\tau_{\delta}$ if and only if $\sha_F = 0$. 

\end{theo}

\noindent
{\bf Proof.} The polynomial $F$ satisfies Condition (C 1), since $F(1) = -1$ and $F(-1) = 1$.

\medskip Let us prove (a). Let $\sigma_1 \in  {\rm Mil}_{1,17}(S)$ and $\sigma_2 \in  {\rm Mil}_{2,2}(C)$ be the restrictions of 
of $\tau_{\zeta}   \in {\rm Mil}_{3,19}(F)$. Since $S$ and $C$ are both irreducible, we have $\sha_S = 0$ and $\sha_C = 0$.
Therefore by Corollary \ref {final coro}, $\Lambda_{1,17}$ has an isometry with characteristic polynomial $S$ and Milnor
index $\sigma_1$ and $\Lambda_{2,2}$ has an isometry with characteristic polynomial $C$ and Milnor
index $\sigma_2$. This implies (a).

\medskip
Let us prove (b). If $\sha_F = 0$, then Corollary \ref {final coro} implies that $\Lambda_{3,19}$ has an isometry with characteristic polynomial $F$ and any Milnor index.

\medskip
Assume that $\sha_F \not = 0$; since $F$ has two irreducible factors, this implies that $\sha_F \simeq {\bf Z}/2{\bf Z}$. 
Recall from \S \ref{local data section} that $\epsilon_{\tau_{\delta}} = \epsilon^{\rm finite} + \epsilon^{\infty}_{\tau_{\delta}}$ and
$\epsilon_{\tau_{\zeta}} = \epsilon^{\rm finite} + \epsilon^{\infty}_{\tau_{\zeta}}$. 

\medskip
By (a) we know that  $\Lambda_{3,19}$ has an isometry with characteristic polynomial $F$ and Milnor
index $\tau_{\zeta}$; this implies that $\epsilon_{\tau_{\zeta}} = 0$. Note that $\epsilon^{\infty}_{\tau_{\delta}} \not = 
\epsilon^{\infty}_{\tau_{\zeta}}$. Therefore $\epsilon_{\tau_{\delta}} \not = 0$, and by Theorem \ref{final} this implies that
$\Lambda_{3,19}$ does not have an isometry with characteristic polynomial $F$ and Milnor
index $\tau_{\delta}$. This completes the proof of (b).

\begin{example} Let $S$ be the Salem polynomial corresponding to the Salem number $\lambda_{18}$. 
This polynomial satisfies the conditions of Proposition 
\ref{resultant} : we have $|{\rm Res}(S,f)| = 1$ for all $f \in \{\Phi_1,\Phi_2,\Phi_3,\Phi_4, \Phi_6 \}$. Therefore
by Proposition \ref{resultant}, if $\Lambda_{3,19}$ has an isometry with characteristic polynomial $S C$ 
and Milnor index $\tau_{\delta}$
for
some product $C$ of cyclotomic polynomials, then we have $C = \Phi_{12}$. 

\medskip
Let $F = S \Phi_{12}$. We have $\sha_F \not = 0$. Indeed, $|{\rm Res}(S,f)| = 169$, and the common factors modulo $13$
of $S$ and $\Phi_{12}$ in ${\bf F}_{13}[X]$ are $X + 6, X + 11 \in {\bf F}_{13}[X]$. These polynomials are not symmetric.
Therefore $\Pi_{S,\Phi_{12}} = \varnothing$, and hence  $\sha_F \simeq {\bf Z}/2{\bf Z}$. Theorem \ref{phi 12} implies
that $\Lambda_{3,19}$ does not have any isometry with characteristic polynomial $S \Phi_{12}$ 
and Milnor index $\tau_{\delta}$.

\medskip
Since this holds for all roots $\delta$ of $S$ with $|\delta| = 1$, the Salem number $\lambda_{18}$ is not
realized by an automorphism of a non-projective $K3$ surface.


\end{example}

\bigskip
\bigskip
Eva Bayer--Fluckiger 

EPFL-FSB-MATH

Station 8

1015 Lausanne, Switzerland

\medskip

eva.bayer@epfl.ch


\begin{thebibliography}{99}








\bibitem[B 15]{B 15} E. Bayer-Fluckiger,
\textit{Isometries of quadratic spaces}, J. Eur. Math. Soc.  \textbf{17} (2015), 1629-1656.

\bibitem[B 20]{B 20} E. Bayer-Fluckiger,
\textit{Isometries of lattices and Hasse principles}, J. Eur. Math. Soc. (to appear),  arXiv:2001.07094.



\bibitem[BT 20]{BT} E. Bayer-Fluckiger, L. Taelman,
\textit{Automorphisms of even unimodular lattices and equivariant Witt groups}, J. Eur. Math. Soc.
\textbf{22} (2020), 3467-3490.

\bibitem[Bo 77]{Bo} D. W. Boyd, \textit{Small Salem numbers}, Duke Math. J. {\bf 44} (1977), 315-328.

 

\bibitem[Br 20]{Br} S. Brandhorst, \textit{On the stable dynamical spectrum of complex surfaces},
Math. Ann. {\bf 377} (2020), 421-434.

\bibitem[BGA 16]{BG} S. Brandhorst, V. Gonzalez-Alonso, \textit{Automorphisms of minimal entropy
on supersingular $K3$ surfaces}, J. London Math. Soc. {\bf 97} (2016), 282-305. 

\bibitem[Ca 14]{Ca} S. Cantat, \textit {Dynamics of automorphisms of compact complex surfaces}, Frontiers in Complex
Dynamics : In celebration of John Milnor's 80th birthday (463-514), Princeton Math. Ser.  {\bf 51}, Princeton Univ. Press,
Princeton, NJ, 2014.

\bibitem[Co 89]{Co} D. Cox, \textit {Primes of the form $x^2 + n y^2$, Fermat, class field theory and complex multiplication},
A Wiley-Interscience Publication, John Wiley \& Sons, Inc., New York, 1989.



\bibitem[E 84]{E} J-H. Evertse, \textit {On equations in $S$-units and the Thue-Mahler equation}, Invent. Math. {\bf 75}
(1984), 561-584.

\bibitem[GM 02]{GM} B. Gross, C. McMullen,
\textit{Automorphisms of even, unimodular lattices and unramified Salem numbers}, J. Algebra {\bf 257} (2002), 265--290. 




\bibitem[H 16]{H} D. Huybrechts, \textit{Lectures on $K3$ surfaces}, Cambridge Studies in Advanced Mathematics {158}, Cambridge
University Press, Cambridge, 2016.


\bibitem [McM 02]{Mc1} C. McMullen, Dynamics on K3 surfaces: Salem numbers and Siegel disks. J. Reine Angew. Math. {\bf 545} (2002), 201-233.

\bibitem[McM 11]{Mc2} C. McMullen, \textit K3 surfaces, entropy and glue. J. Reine Angew. Math. {\bf 658} (2011), 1-25. 

\bibitem[McM 16]{Mc3} C. McMullen, \textit
Automorphisms of projective K3 surfaces with minimum entropy. Invent. Math. {\bf 203} (2016), 179-215.

\bibitem[M 68]{M 68} J. Milnor, \textit {Infinite cyclic coverings}, Topology of Manifolds (J. Hocking, ed.), Prindle, Weber
and Schmidt, Boston (1968), 115–133.

\bibitem[M 69]{M} J. Milnor, \textit {Isometries of inner product spaces}, Invent. Math. {\bf  8}  (1969), 83--97. 




\bibitem[O 10]{O} K. Oguiso, \textit {The third smallest Salem number in automorphisms of $K3$ surfaces}, Algebraic
geometry in East Asia - Seoul 2008, 331-360, Adv. Stud. Pure Math. {\bf 60}, Math. Soc. Japan, Tokyo, 2010.



\bibitem[O'M 73]{OM} O.T. O'Meara, \textit {Introduction to quadratic forms}, reprint of the 1973 edition. Classics in Mathematics. Springer-Verlag, Berlin, 2000.



\bibitem[R 12]{R} P. Reschke, \textit{Salem numbers and automorphisms of complex surfaces}, Math. Res. Lett. {\bf 19}
(2012), 475-482.

\bibitem[R 17]{R2} P. Reschke, \textit{Salem numbers and automorphisms of abelian surfaces}, Osaka J. Math. {\bf 54} (2017),
1-15.



\bibitem[Sch 85]{Sch}  W. Scharlau, \textit {Quadratic and hermitian forms}, Grundlehren der Mathematischen Wissenschaften {\bf 270}, Springer-Verlag, Berlin, 1985.

\bibitem[S 77]{S} J-P. Serre, \textit {Cours d'arithm\'etique}, Presses Universitaires de France, 1977.





\bibitem [Sm 15]{Sm} C. Smyth, \textit {Seventy years of Salem numbers}, Bull. Lond. Math. Soc.
{\bf 47} (2015), 379-395. 

\bibitem [Z 19]{Z} S. Zhao, \textit {Automorphismes loxodromiques de surfaces ab\'eliennes r\'eelles}, Ann. Fac. Sci. Toulouse Math
{\bf 28}, 109-127.



\end{thebibliography}
\end{document}